\documentclass[12pt,a4paper]{article}

\usepackage[left=1cm,right=1cm,top=2cm,bottom=2cm]{geometry}
\usepackage[dvips]{graphics}
\usepackage[dvips]{epsfig}
\usepackage{color}
\usepackage{hyperref}
\usepackage{tikz,authblk}

\newtheorem{lemma}{Lemma}
\newtheorem{thm}{Theorem}
\newtheorem{cor}[lemma]{Corollary}
\newtheorem{prop}[lemma]{Proposition}
\newtheorem{defn}{Definition}
\newtheorem{rem}{Remark}
\newtheorem{conj}{Conjecture}

\newcommand{\dimo}{\noindent \emph{Proof. }}
\newcommand{\qed}{\\ \rightline{$\Box$ \ \ \ \ \ \ \ \ \ \ \ \ \ \ \ }\\}
\newcommand{\e}{\varepsilon}

\usepackage{latexsym}
\usepackage{amsfonts}
\usepackage{amssymb}
\usepackage{amsmath}
\usepackage[english]{babel}

\begin{document}

 \title{CLASSIFYING COMPACT 4-MANIFOLDS VIA GENERALIZED REGULAR GENUS AND G-DEGREE}

 \renewcommand{\Authfont}{\scshape\small}
 \renewcommand{\Affilfont}{\itshape\small}
 \renewcommand{\Authand}{ and }

\author[1] {Maria Rita Casali}
\author[2] {Paola Cristofori}

 \affil[1] {Department of Physics, Mathematics and Computer Science, University of Modena and Reggio Emilia \ \ \ \ \ \ \ \ \ \ \   Via Campi 213 B, I-41125 Modena (Italy), casali@unimore.it}

\affil[2] {Department of Physics, Mathematics and Computer Science, University of Modena and Reggio Emilia,  \ \ \ \ \ \ \ \ \ \   Via Campi 213 B, I-41125 Modena (Italy),  paola.cristofori@unimore.it}

\maketitle

\begin{abstract}
{\it $(d+1)$-colored graphs}, i.e. edge-colored graphs that are $(d+1)$-regular, have already been proved to be a useful representation tool for compact PL $d$-manifolds, thus extending the theory (known as {\it crystallization theory}) originally developed for the closed case. In this context, combinatorially defined PL invariants  play a relevant role. 
The present paper focuses in particular on {\it generalized regular genus} and {\it G-degree}: the first one extending to higher dimension the classical notion of Heegaard genus for 3-manifolds, the second one arising, within theoretical physics, from the theory of random tensors as an approach to quantum gravity in dimension greater than two. \\
We establish several general results concerning the two invariants, in relation with invariants of the boundary and with the rank of the fundamental group, as well as their behaviour with respect to connected sums.
We also compute both generalized regular genus and G-degree for interesting classes of compact $d$-manifolds, such as handlebodies, products of closed manifolds by the interval and $\mathbb D^2$-bundles over $\mathbb S^2.$ \\
The main results of the paper concern dimension 4, where we 
obtain the classification of all compact PL manifolds with generalized regular genus at most one, and of all compact PL manifolds with G-degree at most 18; moreover, in case of empty or connected boundary, the classifications are extended to generalized regular genus two and to G-degree 24. 
\end{abstract}



\section{Introduction}\label{intro}

In the PL $d$-dimensional setting ($d \ge 3$), both the invariants {\it generalized regular genus} and {\it G-degree} have been recently introduced, making use of the possibility of representing all compact PL $d$-manifolds by means of regular $(d+1)$-colored graphs (i.e. graphs whose vertices have degree $d+1$, and so that the $d+1$ edges adjacent to each vertex are injectively colored by the colors $\{0,1, \dots, d \}$): see  \cite{Casali-Cristofori-Grasselli} and \cite{Casali-Cristofori-Dartois-Grasselli} respectively, or the following Section  \ref{prelim}. 

This kind of representation for compact PL manifolds and the study of the above invariants has been deeply motivated by the strong connections - accurately described in \cite{Casali-Cristofori-Dartois-Grasselli} - between 
random tensor models and the so called {\it crystallization theory}, which is a useful combinatorial tool for the topological and geometrical study of PL manifolds of arbitrary dimension (assumed to be closed, in the ``classical" version of the theory)  via edge-colored graphs.  

Without going into details, we only recall that colored tensor models  were introduced within theoretical physics as a random geometry approach to quantum gravity, following the successes of random matrix models (\cite{DiFrancesco-al})  
and, in this context, the coefficients of the $1/N$ expansion of the correlation functions   
in dimension $d$ are generating functions of regular bipartite $(d+1)$-colored graphs; 
moreover, the quantity driving the $1/N$ expansion is the {\it Gurau degree}, whose definition involves  the genera of the surfaces where the considered graphs regularly embed (see Definition \ref{reg_emb} for details), 
exactly as  the strictly related notions of generalized regular genus and G-degree (\cite{[uncoloring]}, \cite{[Bonzom-Gurau-Riello-Rivasseau]}, \cite{[Gurau-Schaeffer 2013]}).    
Hence, any result obtained about generalized regular genus and/or G-degree is not only an achievement in the comprehension and possible classification of manifolds in the PL category, 
but may also bring insights in the interactions between geometry and physics.\footnote{Note that, while the so-called colored tensor model is usually formulated in terms of complex tensors, a version involving real tensor fields, known as the Gurau-Witten model, has 
gained considerable interest in the study of quantum mechanical models with random interactions and their holographic properties (\cite{Gurau2017}, \cite{Witten}). 
Its $1/N$ expansion involves both bipartite and non-bipartite graphs, and hence the investigation concerns both orientable and non-orientable manifolds.}

Other types of tensor models whose Feynman graphs include or are a particular kind of colored graphs have also been introduced (for instance the O(N)$^3$-invariant, the multi-orientable or the colored SYK models; see \cite{[Bonzom-Lionni-Tanasa]}, \cite{Carrozza-Tanasa}, \cite{Gurau2019}, \cite{Tanasa-multiorientable} and references therein)
and involve the definition of suitable ``degrees'' on colored graphs whose topological significance with respect to the dual colored triangulations can be investigated,
 as done for example in \cite{Fusy-Lionni-Tanasa}.
\medskip

As far as the $d$-dimensional setting is concerned, the present paper proves in Section \ref{sec_general-properties} several general properties of generalized regular genus and G-degree, relating them to the analogue invariants for the boundary manifold (Proposition \ref{genus vs boundary genus}), or for the summands of a connected sum decomposition (Proposition   \ref{connected sums}), and establishing an inequality between generalized regular genus and the rank of the fundamental group, in case of manifolds with empty or connected boundary (Proposition   \ref{genus vs rank}).  
Moreover, in Section  \ref{sec_hadlebodies and xI},  standard graphs representing some interesting classes of PL $d$-manifolds are obtained, yielding the computation of the associated invariants: see Proposition \ref{prop.handlebodies} concerning handlebodies, and Proposition \ref{MxI} concerning the product between a closed $d$-manifold and the interval.  
A similar approach is then performed in Section \ref{section_bundles}, in the $4$-dimensional setting, as regards the $\mathbb D^2$-bundles over $\mathbb S^2$. 

The focus of the paper is, indeed, on dimension $d=4$, where crystallization theory already yielded classification results in the closed case both with respect to regular genus and with 
respect to gem-complexity (see \cite{Casali-Cristofori ElecJComb 2015}, \cite{Casali-Cristofori-Gagliardi Complutense 2015} and their references).
\\
On the whole, experimental approaches to PL classification of 4-manifolds via combinatorial descriptions and associated PL-invariants are recent and still developing: 
as interesting examples we recall the use of Turaev's shadows and shadow complexity (see for instance \cite{Koda-Martelli-Naoe}, \cite{Martelli} for the closed case, and
\cite{Naoe-Proc.AMS} for the case of compact acyclic 4-manifolds) and that of trisections and trisection genus (see \cite{Meier}, \cite{Meier-Zupan}, 
\cite{Spreer-Tillmann})\footnote{Generalizing \cite{Spreer-Tillmann}, a recent study also performs an approach to trisection genus via crystallizations, giving rise 
to the notion of {\it gem-induced trisection genus}, which applies also to manifolds with connected boundary: see \cite{[Casali-Cristofori_trisection]}.}.

In the 4-dimensional setting, the present paper applies the general combinatorial properties of graphs representing compact $d$-manifold (obtained in Section \ref{Sec-comb_properties}), 
together with classical methods of crystallization theory and recent achievements about Dehn surgery, in order to yield classifying results for compact PL $4$-manifolds $M^4$ with respect to both their generalized regular genus $\bar{ \mathcal G}(M^4)$ and their G-degree $\mathcal D_G(M^4)$. 

In particular, we prove in Section \ref{section.4dim} the following  statements (where $\mathbb S^1 \times \mathbb S^3$ and $\mathbb S^1 \tilde \times \mathbb S^3$ denote the orientable and non-orientable sphere bundle over $\mathbb S^1$, ${\mathbb Y}^4_m$ and $ \tilde{ \mathbb Y}^4_m$ denote the orientable and non-orientable $4$-handlebody of genus $m$, 
$\xi_c$ denotes the $\mathbb D^2$-bundle over $\mathbb S^2$ with Euler class $c$, while $M^4(K,d)$ denotes the compact PL $4$-manifold obtained from the $4$-disk by adding a $2$-handle according to the framed knot $(K,d)$, whose boundary is the $3$-manifold $M^3(K,d)$ obtained from the $3$-sphere by Dehn surgery on the same framed knot):  

\begin{thm} \label{Thm.gen-genus_one&two}
Let $M^4$ be a compact PL 4-manifold with no spherical boundary components. Then:  
\begin{itemize} 
\item[a.] \ $\bar{ \mathcal G}(M^4)=0$ \ if and only if  \ $M^4 \cong\mathbb S^4.$ 
\item[b.] \ $\bar{ \mathcal G}(M^4)=1$ \ if and only if  \ 
\begin{itemize}  
\item[$\bullet$] either \ $M^4 \in \{\mathbb S^1 \times \mathbb S^3, \   \mathbb S^1 \tilde \times \mathbb S^3\}$,  
\item[$\bullet$]  or \  $M^4 \in \{{\mathbb Y}^4_1, \tilde{  \mathbb Y}^4_1\}$,   
\item[$\bullet$]  or \  $M^4 \cong \bar M \times I$, \ where $\bar M$ is a genus one closed 3-manifold; 
\end{itemize}
\item[c.] If $M^4$ has empty or connected boundary and \ $\bar{ \mathcal G}(M^4)=2,$ then: 
\begin{itemize}  
\item[$\bullet$] either \ $M^4 \in \{\#_2 (\mathbb S^1 \times \mathbb S^3), \  \#_2(\mathbb S^1 \tilde \times \mathbb S^3), \  \mathbb{CP}^2 \},$   
\item[$\bullet$] or  \ $M^4 \in \{{\mathbb Y}^4_2, \ \tilde{ \mathbb Y}^4_2, 
 (\mathbb S^1 \times \mathbb S^3)\# {\mathbb Y}^4_1, \  (\mathbb S^1 \tilde \times \mathbb S^3)\# {\mathbb Y}^4_1, \  (\mathbb S^1 \tilde \times \mathbb S^3)\# \tilde{\mathbb Y}^4_1,  
\  \mathbb S^2 \times \mathbb D^2, \ \xi_c\}$, 
\item[$\bullet$]  or \ $M^4 \cong M^4(K,d),$ \ $(K,d)$ being a non-trivial framed knot such that \ $M^3(K,d) =  L(\alpha, \beta)$ with $\alpha \ge 3.$ \ 
\end{itemize}
\end{itemize}
\end{thm}

\begin{thm} \label{thm.classif_G-degree}
Let $M^4$ be a compact PL $4$-manifold  with no spherical boundary components. Then:  
\begin{itemize}
 \item [a.] \ $\mathcal D_G(M^4)  =  0$ \ if and only if  \ $M^4 \cong\mathbb S^4;$
 \item [b.] \ $\mathcal D_G(M^4) = 12$ \ if and only if  \ 
\begin{itemize}  
\item[$\bullet$]  either \  $M^4 \in \{\mathbb S^1 \times \mathbb S^3, \ \mathbb S^1 \tilde \times \mathbb S^3\}$, 
 \item[$\bullet$]   or  \ $M^4 \in \{   \mathbb Y^4_1, \ \tilde{\mathbb Y}^4_1\};$
\end{itemize} 
 \item [c.] \ $\mathcal D_G(M^4)  = 18$ \ if and only if  \ $ M^4 \in \{L(2,1) \times I, \  (\mathbb S^1 \times \mathbb S^2) \times I,   \  (\mathbb S^1 \tilde \times \mathbb S^2) \times I \}$. 
 \end{itemize}
\noindent No other compact PL $4$-manifold (with no spherical boundary components) exists with $\mathcal D_G(M^4)\leq 23.$     
\vskip 5pt
\noindent 
Moreover, if $M^4$ has empty or connected boundary, then: 
\begin{itemize}
 \item [d.] \ $ \mathcal D_G(M^4) = 24 $ \ if and only if  \ 
\begin{itemize}  
\item[$\bullet$]  either \ $M^4 \in \{\#_2(\mathbb S^1 \times \mathbb S^3), \  \#_2(\mathbb S^1 \tilde \times \mathbb S^3), \  \mathbb{CP}^2 \}$, 
\item[$\bullet$]  or \ $M^4\in \{{\mathbb Y}^4_2, \ \tilde{ \mathbb Y}^4_2,  
 (\mathbb S^1 \times \mathbb S^3)\# {\mathbb Y}^4_1, \  (\mathbb S^1 \tilde \times \mathbb S^3)\# {\mathbb Y}^4_1, \  (\mathbb S^1 \tilde \times \mathbb S^3)\# \tilde{\mathbb Y}^4_1,  
\  \mathbb S^2 \times \mathbb D^2, \ \xi_2\}.$ 
\end{itemize} 
 \end{itemize}
 \end{thm}
 
A direct consequence of Theorem \ref{thm.classif_G-degree} is the identification of all  compact orientable PL $4$-manifolds (resp. compact orientable  PL $4$-manifolds with empty or connected boundary), 
represented by regular graphs involved in the first four (resp. five) most significant non-null terms of the $1/N$ expansion
of the correlation functions (\cite{[uncoloring]}, \cite{Casali-Grasselli 2017}): see Remark \ref{1/N-expansion}. 
 
Further results of the present paper are also the characterization of $4$-dimensional handlebodies as the only PL $4$-manifolds with connected (non empty) boundary whose 
generalized regular genus equals that of their boundary (Theorem   \ref{genus vs boundary genus (n=4)}),   
while the equality between generalized regular genus and the rank of the fundamental 
group characterizes $\#_{\rho}(\mathbb S^1 \times \mathbb S^3)$ and $\#_{\rho}(\mathbb S^3 \widetilde{\times} \mathbb S^1)$ in the closed case,  
$\#_{\alpha} (\mathbb S^1 \times \mathbb S^3)\# {\mathbb Y}^4_{\beta},$ \  $\#_{\alpha} (\mathbb S^1 \tilde \times \mathbb S^3)\# {\mathbb Y}^4_{\beta},$ \  $\#_{\alpha} (\mathbb S^1 \tilde \times \mathbb S^3)\# \tilde{\mathbb Y}^4_{\beta}$,   
in the connected boundary case (Theorem \ref{genus vs rank (n=4)}). 

Note that, as a consequence of the above results, all $\mathbb D^2$-bundles over $\mathbb S^2$ turn out to have generalized regular genus 2, thus proving that  generalized regular genus is not finite-to-one in dimension four, and that, in the $4$-dimensional case of non-empty  boundary,  
the equality between generalized regular genus and the ``classical"  invariant {\it regular genus} does not hold, even if the boundary is assumed to be connected (Corollary \ref{cor.no-equality (n=4)}).

\section{Preliminaries}  \label{prelim}

In the present section we will briefly review some basic notions of the so called {\it crystallization theory}, which is a representation tool for general piecewise linear (PL) compact manifolds, without assumptions about dimension, connectedness, orientability or boundary properties  (see the ``classical" survey paper \cite{Ferri-Gagliardi-Grasselli}, or the more recent one \cite{Casali-Cristofori-Gagliardi Complutense 2015}, concerning the $4$-dimensional case).

From now on, unless otherwise stated, all spaces and maps will be considered in the PL category, and all manifolds will be assumed to be compact and connected.

\begin{defn} \label{$n+1$-colored graph}
{\em A {\it $(d+1)$-colored graph}  ($d \ge 2$) is a pair $(\Gamma,\gamma)$, where $\Gamma=(V(\Gamma), E(\Gamma))$ is a multigraph (i.e. multiple edges are allowed, but loops are forbidden) 
which is regular of degree  $d+1$, and $\gamma$ is an {\it edge-coloration}, that is a map  $\gamma: E(\Gamma) \rightarrow \Delta_d=\{0,\ldots,d\}$ which is injective on 
adjacent edges.}\end{defn}

In the following, for sake of concision, when the coloration is clearly understood, we will denote colored graphs simply by $\Gamma$. 

\smallskip

For every $\{c_1, \dots, c_h\} \subseteq\Delta_d$ let $\Gamma_{\{c_1, \dots, c_h\}}$  be the subgraph obtained from $(\Gamma, \gamma)$ by deleting all the edges that are not colored by the elements of $\{c_1, \dots, c_h\}$. 
In this setting, the complementary set of $\{c\}$ (resp. $\{c_1,\dots,c_h\}$)  in $\Delta_d$ will be denoted by $\hat c$ (resp. $\hat c_1\cdots\hat c_h$). 
The connected components of $\Gamma_{\{c_1, \dots, c_h\}}$ are called {\it $\{c_1, \dots, c_h\}$-residues} or {\it $h$-residues} of $\Gamma$; their number is denoted by $g_{\{c_1, \dots, c_h\}}$ (or, for short, by $g_{c_1,c_2}$, $g_{c_1,c_2,c_3}$ and $g_{\hat c}$ if $h=2,$ $h=3$ and $h = d$ respectively). 

\medskip 

\noindent A $d$-dimensional pseudocomplex $K(\Gamma)$ can be associated to a $(d+1)$-colored graph $\Gamma$: 
\begin{itemize}
\item take a $d$-simplex for each vertex of $\Gamma$ and label its vertices by the elements of $\Delta_d$;  
\item if two vertices of $\Gamma$ are $c$-adjacent ($c\in\Delta_d$), glue the corresponding $d$-simplices  along their $(d-1)$-dimensional faces opposite to the $c$-labeled vertices, so that equally labeled vertices are identified.
\end{itemize}

\smallskip

In general $|K(\Gamma)|$ is a {\it $d$-pseudomanifold} and $\Gamma$ is said to {\it represent} it. 

\medskip

Note that, by construction, $K(\Gamma)$ is endowed with a vertex-labeling by $\Delta_d$ that is injective on any simplex. Moreover, $\Gamma$ turns out to be the 1-skeleton of the dual complex of $K(\Gamma)$.
The duality establishes a bijection between the $\{c_1, \dots, c_h\}$-residues of  $\Gamma$ and the $(d-h)$-simplices of $K(\Gamma)$ whose vertices are labeled by   $\Delta_d - \{c_1, \dots, c_h\}$. 

Given a pseudocomplex $K$ and an $h$-simplex $\sigma^h$ of $K$, the {\it disjoint star} of $\sigma^h$ in $K$ is the pseudocomplex obtained by taking all $d$-simplices of $K$ having $\sigma^h$ as a face
and identifying only their faces that do not contain $\sigma^h.$ The {\it disjoint link}, $lkd(\sigma^h,K)$, of $\sigma^h$ in $K$ is the subcomplex of the disjoint star formed
by those simplices that do not intersect $\sigma^h.$

\noindent 

In particular, given a $(d+1)$-colored graph $\Gamma$, each connected component of $\Gamma_{\hat c}$ ($c\in\Delta_d$) is a $d$-colored graph representing the disjoint link of a $c$-labeled vertex of $K(\Gamma)$, 
that is also (PL) homeomorphic to the link of this vertex in the first barycentric subdivision of $K(\Gamma).$

\begin{defn} {\em A {\it singular (PL) $d$-manifold} is a closed connected $d$-dimensional polyhedron admitting a simplicial triangulation where the links of vertices are closed connected $(d-1)$-manifolds, 
while  the links of all $h$-simplices of the triangulation with $h > 0$ are (PL) $(d-h-1)$-spheres. Vertices whose links are not PL  $(d-1)$-spheres are called {\it singular}.  
\\ \noindent Note that, in case of polyhedra arising from colored graphs, the condition about links of vertices obviously implies the one about links of $h$-simplices, with $h> 0.$}
\end{defn}

\noindent Therefore: 
\begin{itemize}
\item  $|K(\Gamma)|$ is a singular $d$-manifold iff, for each color $c\in\Delta_d$, all $\hat c$-residues of $\Gamma$ represent closed connected $(d-1)$-manifolds. 
\end{itemize}

\noindent In particular: 
\begin{itemize}
\item $|K(\Gamma)|$ is a closed $d$-manifold iff, for each color $c\in\Delta_d$, all $\hat c$-residues of $\Gamma$ represent the $(d-1)$-sphere.
\end{itemize}

\begin{rem} \label{correspondence-sing-boundary} {\em If $N$ is a singular $d$-manifold, then a compact $d$-manifold $\check N$ is easily obtained by deleting small open neighbourhoods of its singular vertices.
Obviously $N=\check N$ iff $N$ is a closed manifold; otherwise, $\check N$ has non-empty boundary (without spherical components).
Conversely, given a compact $d$-manifold $M$, a singular $d$-manifold $\widehat M$ can be constructed by capping off each component of $\partial M$ by a cone over it.

Note that, by restricting ourselves to the class of compact $d$-manifolds with no spherical boundary components,  the above correspondence is bijective and so singular $d$-manifolds and compact $d$-manifolds of this class can be associated  to each other in a well-defined way.

For this reason, throughout the present work, we will restrict our attention to compact manifolds without spherical boundary components. Obviously, in this wider context, closed $d$-manifolds are characterized by $M= \widehat M.$}
\end{rem}

In virtue of the bijection described in Remark \ref{correspondence-sing-boundary}, a $(d+1)$-colored graph $\Gamma$ is said to {\it represent}
a compact $d$-manifold $M$ with no spherical boundary components if and only if  it represents the associated singular manifold $\widehat M$.

\medskip
The following theorem extends to the boundary case a well-known result - originally stated in \cite{Pezzana}  - founding the combinatorial representation theory for closed manifolds of arbitrary dimension via colored graphs.

\begin{prop}[\cite{Casali-Cristofori-Grasselli}]\ \label{Theorem_gem}  
Any compact orientable (resp. non orientable) $d$-manifold with no spherical boundary components admits a bipartite (resp. non-bipartite) $(d+1)$-colored graph representing it.
\end{prop}

If $\Gamma$ represents  a compact $d$-manifold,  a $d$-residue of $\Gamma$ will be called {\it ordinary} if it represents $\mathbb S^{d-1}$, {\it singular} otherwise. Similarly, a color $c$ will be called {\it singular} if at least one of the $\hat c$-residues of $\Gamma$ is singular.

\medskip

The  existence of a particular type of embedding of colored graphs into surfaces, is the key result in order to define the important notion of  regular genus.

\begin{prop}[\cite{Gagliardi 1981}]\label{reg_emb}   
Let $(\Gamma,\gamma)$ be a bipartite (resp. non-bipartite) $(d+1)$-colored graph of order $2p$. Then for each cyclic permutation $\varepsilon = (\varepsilon_0,\ldots,\varepsilon_d)$ of $\Delta_d$, up to inverse, there exists a cellular embedding, called {\rm regular}
\footnote{Regular embeddings are called \emph{Jackets} in the tensor models context.}, of $(\Gamma,\gamma)$ into an orientable (resp. non-orientable) closed surface $F_{\varepsilon}(\Gamma)$
whose regions are bounded by the images of the $\{\varepsilon_j,\varepsilon_{j+1}\}$-colored cycles, for each $j \in \mathbb Z_{d+1}$.
Moreover, the genus (resp. half the genus)  $\rho_{\varepsilon} (\Gamma)$ of $F_{\varepsilon}(\Gamma)$ satisfies

\begin{equation*}
2 - 2\rho_\varepsilon(\Gamma)= \sum_{j\in \mathbb{Z}_{d+1}} g_{\varepsilon_j\varepsilon_{j+1}} + (1-d)p.
\end{equation*}

No regular embeddings of $(\Gamma,\gamma)$ exist into non-orientable (resp. orientable) surfaces.
\end{prop}

The \emph{Gurau degree} (often called {\it degree} in the tensor models literature, see \cite{Gurau-book}) and the {\it regular genus} of a colored graph are defined in terms of the embeddings of Proposition \ref{reg_emb}. 

\begin{defn} \label{Gurau-degree}
{\em Let $(\Gamma,\gamma)$ be a $(d+1)$-colored graph.
If $\{\varepsilon^{(1)}, \varepsilon^{(2)}, \dots , \varepsilon^{(\frac {d!} 2)}\}$ is the set of all cyclic permutations of $\Delta_d$ (up to inverse), $ \rho_{\varepsilon^{(i)}}(\Gamma)$ \ ($i=1, \dots , \frac {d!} 2$) is called the {\it regular genus of $\Gamma$ with respect to the permutation $\varepsilon^{(i)}$}. 
Then, the \emph{Gurau degree} (or \emph{G-degree} for short) of $\Gamma$, denoted by  $\omega_{G}(\Gamma)$, is defined as}
\begin{equation*}
 \omega_{G}(\Gamma) \ = \ \sum_{i=1}^{\frac {d!} 2} \rho_{\varepsilon^{(i)}}(\Gamma)
\end{equation*}
{\em and the {\it regular genus} of $\Gamma$, denoted by $\rho(\Gamma)$, is defined as}
\begin{equation*}
 \rho(\Gamma) \ = \ \min\, \{\rho_{\varepsilon^{(i)}}(\Gamma)\ /\ i=1,\ldots,\frac {d!} 2\}.
\end{equation*}
\end{defn}

Note that, in dimension $2$, any bipartite (resp. non-bipartite) $3$-colored graph $(\Gamma,\gamma)$ represents an orientable (resp. non-orientable) surface $|K(\Gamma)|$ and $\rho(\Gamma)= \omega_G(\Gamma)$ is exactly the genus (resp. half the genus) of $|K(\Gamma)|.$
On the other hand, for $d\geq 3$, the G-degree of any $(d+1)$-colored graph (resp. the regular genus of any $(d+1)$-colored graph representing a closed PL $d$-manifold) is proved to be a non-negative {\it integer}, both in the bipartite and non-bipartite case:  see \cite[Proposition 7]{Casali-Cristofori-Dartois-Grasselli} (resp. \cite[Proposition A]{Chiavacci-Pareschi}).

\medskip

As a consequence of the definition of regular genus of a colored graph and of Proposition \ref{Theorem_gem}, two PL invariants for compact $d$-manifolds can be defined:  

\begin{defn}\label{def_gen_degree} {\em Let $M$ be a compact (PL) $d$-manifold ($d\geq 2$).
The {\it generalized regular genus} of $M$ is defined as}
\begin{equation*}
\overline{\mathcal G}(M)=\min \{\rho(\Gamma)\ | \ (\Gamma,\gamma)\mbox{ represents} \ M\}.
\end{equation*}
{\em and the {\it Gurau degree} (or {\it G-degree}) of $M$ is defined as}
\begin{equation*}
\mathcal D_G(M)=\min \{\omega_G(\Gamma)\ | \ (\Gamma,\gamma)\mbox{ represents} \ M\}.
\end{equation*}
\end{defn}

\medskip

For any $(d+1)$-colored graph $\Gamma$,  the following inequality obviously holds: \ 
$ \mathcal \omega_G(\Gamma) \ \ge \ \frac{d!}2 \cdot \rho(\Gamma).$ 
Hence, for any compact $d$-manifold $M$: 
$$ \mathcal D_G(M) \ \ge \ \frac{d!}2 \cdot \overline{\mathcal G}(M).$$

\begin{rem} \label{rem_regular genus} {\em Note that, 
in case $M$ is a closed PL $d$-manifold, the generalized regular genus coincides by definition with the PL invariant {\it regular genus} (see Section 4),   extending to higher dimension the Heegaard genus of a $3$-manifold (\cite{Cristofori-Gagliardi-Grasselli}). 
Regular genus zero  succeeds in characterizing spheres in arbitrary dimension,\footnote{In Proposition \ref{regular genus vs generalized regular genus} we will prove that the same characterization holds for generalized regular genus, too, within the wider class of compact $d$-manifolds.}  and many  
classifying results via regular genus have been obtained, especially in dimension $4$ and $5$ (see \cite{Basak-Casali 2016}, \cite{Casali_Forum2003}, \cite{Casali-Cristofori-Gagliardi Complutense 2015}, \cite{Casali-Gagliardi ProcAMS} and their references). 
Also G-degree zero characterizes spheres in arbitrary dimension, and some classifying results via this invariant have recently been obtained for compact 3-manifolds and for closed PL 4-manifolds: see \cite{Casali-Cristofori-Dartois-Grasselli} and \cite{Casali-Cristofori-Grasselli}.}
\end{rem}  

Finally, we recall that, within crystallization theory, a finite set of combinatorial moves have been defined, which translate the homeomorphism problem   of the represented polyhedra. 

\begin{defn}\label{def_dipole} {\em An {\it $r$-dipole ($1\le r\le d$) of colors $c_1,\ldots,c_r$} of a $(d+1)$-colored graph $(\Gamma,\gamma)$ is a subgraph of $\Gamma$ consisting in two vertices joined by $r$ edges, colored by $c_1,\ldots,c_r$, such that its vertices belong to different connected components of $\Gamma_{\hat c_1\ldots\hat c_r}$.

The {\it elimination} of an $r$-dipole in $\Gamma$ can be performed by deleting the subgraph and welding the remaining hanging edges according to their colors; in this way another $(d+1)$-colored graph $(\Gamma^\prime,\gamma^\prime)$ is obtained. The inverse operation is called the {\it addition} of the dipole to $\Gamma^\prime.$

The dipole is called {\it proper} if $|K(\Gamma)|$ and $|K(\Gamma^\prime)|$ are PL homeomorphic. It is known that this happens when at least one of the two connected components of $\Gamma_{\hat c_1\ldots \hat c_r}$ 
intersecting the dipole represents a $(d-r)$-sphere (\cite[Proposition 5.3]{Gagliardi 1987}).\footnote{Note that, if $|K(\Gamma)|$ is a singular $d$-manifold, all $r$-dipoles are proper, for $1 < r \le d$.}}
\end{defn}

\begin{rem}\label{rem_dipoli} {\em Neither the G-degree nor the regular genus of a $(d+1)$-colored graph are affected by elimination of $1$-dipoles. 
Therefore, from any $(d+1)$-colored graph $\Gamma$ representing a compact PL $d$-manifold $M$ with empty or connected boundary, by eliminating (proper) 1-dipoles, a $(d+1)$-colored graph can be obtained, still representing $M$, with the same G-degree and regular genus as $\Gamma$ and having only one ${\hat\imath}$-residue for each $i\in\Delta_d.$
Such a $(d+1)$-colored graph is said to be a {\it crystallization} of $M$.  }
\end{rem}  

\section{Combinatorial properties of graphs representing singular $d$-manifolds} 
\label{Sec-comb_properties}

In \cite{Casali - Forum Math 1992}, \cite{Casali Canadian} and  \cite{Casali-Gagliardi ProcAMS},  interesting combinatorial formulae have been obtained, regarding both regular edge-colored graphs representing closed $d$-manifolds and  edge-colored graphs with boundary (see 
\cite{Ferri-Gagliardi-Grasselli}, or the next Section \ref{sec_general-properties}) representing $d$-manifolds with non-empty boundary. 
Here, we will generalize them to regular edge-colored graphs representing (via singular $d$-manifolds) all compact (PL) $d$-manifolds. 

\smallskip

In the following, let $(\Gamma,\gamma)$ be a (possibly disconnected)   $(d+1)$-colored graph representing a (possibly disconnected) singular $d$-manifold $N^d$.
If $\Gamma$ (resp. $\Gamma_{\mathcal B}$, with $\mathcal B \subset \Delta_d$) has $g \ge 1$ (resp. $g_{\mathcal B} \ge 1$) connected components 
$\Gamma^1, \Gamma^2, \dots, \Gamma^g$ (resp. $H^1, H^2, \dots, H^{g_{\mathcal B}}$), for each permutation 
$\varepsilon = (\varepsilon_0, \varepsilon_1, \dots , \varepsilon_{d-1}, \varepsilon_d)$ of $\Delta_d$ we define  
$$\rho_\varepsilon(\Gamma)= \sum_{i=1}^g \rho_\varepsilon(\Gamma^i)$$ (resp. $$\rho_{\varepsilon}(\Gamma_{\mathcal B})= \
\sum_{i=1}^{g_{\mathcal B}} \rho_{\varepsilon}(H^i),$$
where by $\rho_{\varepsilon}(H^i)$ we denote the regular genus of $H^i$ with respect to the permutation induced by $\varepsilon$ on the subset $\mathcal B$ of $\Delta_d)$. 

\begin{prop} \label{combinatorial-properties}
If $(\Gamma,\gamma)$ is a  $(d+1)$-colored graph with $g \ge 1$ connected components, representing a (possibly disconnected) singular $d$-manifold $N^d$, then 
\begin{itemize}
\item 
if $\#\mathcal B=m$ and $m \le d-1$, $m$ odd: 
\begin{equation} \label{residue_odd}
2 g_{\mathcal B} = (2-m) p + \sum_{s=2}^{m-1} (-1)^s \sum_{i_1, i_2, \dots, i_s \in \mathcal B} g_{i_1, i_2, \dots, i_s}
\end{equation}
\item 
if $\#\mathcal B=m$ and $m \le d-1$, $m$ even: 
\begin{equation} \label{residue_even}
 0 = (2-m) p + \sum_{s=2}^{m-1} (-1)^s \sum_{i_1, i_2, \dots, i_s \in \mathcal B} g_{i_1, i_2, \dots, i_s}
 \end{equation}
\item 
if $\mathcal B=\Delta_d-\{i\},$ with $i$ non-singular color and $d$ odd: 
\begin{equation} \label{residue_d_odd}
2 g_{\hat\imath} = (2-d) p + \sum_{s=2}^{d-1} (-1)^s \sum_{i_1, i_2, \dots, i_s \in \Delta_d - \{i\}} g_{i_1, i_2, \dots, i_s}
\end{equation}
\item 
if $\mathcal B=\Delta_d-\{i\},$ with $i$ non-singular color and $d$ even: 
\begin{equation} \label{residue_d_even}
 0 = (2-d) p + \sum_{s=2}^{d-1} (-1)^s \sum_{i_1, i_2, \dots, i_s \in \Delta_d - \{i\}} g_{i_1, i_2, \dots, i_s}
 \end{equation} 
\end{itemize}
Moreover: 
\begin{equation} \label{3-residuesversus2-residues}
g_{\varepsilon_{i-1},\varepsilon_{i+1}} =  g_{\varepsilon_{i-1}, \varepsilon_{i}, \varepsilon_{i+1}}  + (\rho_\varepsilon(\Gamma) - \rho_{\varepsilon}(\Gamma_{\widehat {\varepsilon_i}})) - (g - g_{\widehat{\varepsilon_i}}) \ \ \ \ \forall i \in \Delta_d; 
\end{equation}
\begin{eqnarray} \label{numerospigoli}
&&g_{\widehat{\varepsilon_{i}},\widehat{\varepsilon_{j}}} =  (g_{\widehat{\varepsilon_{i}}} + g_{\widehat{\varepsilon_{j}}} -g ) +  \rho_\varepsilon(\Gamma) - \rho_{\varepsilon}(\Gamma_{\widehat {\varepsilon_i}}) - \rho_{\varepsilon
}(\Gamma_{\widehat {\varepsilon_j}}) + \rho_{\varepsilon
}(\Gamma_{\widehat {\varepsilon_j}, \widehat {\varepsilon_i}})\\
&&\forall i,j \text{\ non consecutive in \ } \Delta_d;\nonumber
\end{eqnarray}
\begin{eqnarray} \label{2-residuesversuspairsof3-residues}
g_{\varepsilon_{i-1},\varepsilon_{i+1}} =&&
g_{\varepsilon_{i-1}, \varepsilon_{i}, \varepsilon_{i+1}}  + g_{\varepsilon_{i-1}, \varepsilon_{i+1}, \varepsilon_{r}} + \rho_{\varepsilon}(\Gamma_{\varepsilon_{i-1}, 
\varepsilon_{i}, \varepsilon_{i+1}, \varepsilon_{r}}) + \\ &&- g_{\varepsilon_{i-1}, \varepsilon_{i}, \varepsilon_{i+1}, \varepsilon_{r}}
\qquad \forall r \in \Delta_d - \{i-1,i,i+1\}.\nonumber
\end{eqnarray}
\end{prop}

\dimo
By definition of generalized regular genus       
with respect to the permutation $\varepsilon$:
\begin{equation} \label{genereGamma}
\sum_{j\in \mathbb{Z}_{d+1}} g_{\varepsilon_j,\varepsilon_{j+1}} + (1-d)p = 2 g - 2\rho_\varepsilon(\Gamma)
\end{equation}

By applying the same relation to the (possibly disconnected) subgraph $\Gamma_{\widehat{\varepsilon_i}}$ ($i \in \Delta_d$), we have:
\begin{equation} \label{genereGamma_c}
\sum_{j\in \mathbb{Z}_{d+1}-\{i-1, i\}} g_{\varepsilon_j,\varepsilon_{j+1}} + g_{\varepsilon_{i-1},\varepsilon_{i+1}} + (2-d)p = 2 g_{\widehat{\varepsilon_i}} - 2\rho_{\varepsilon}(\Gamma_{\widehat {\varepsilon_i}})
\end{equation}

In order to prove relations \eqref{residue_odd} and \eqref{residue_even},  recall that each connected component of  $\Gamma_{\mathcal B}$ represents the disjoint link of a $(d-m)$-simplex in the singular $d$-manifold $|K(\Gamma)|=N^d$, which - under the hypothesis $m \le d-1$ - is homeomorphic to the $(m-1)$-sphere. Hence, its Euler characteristic equals $2$ if $m$ is odd and $0$ if $m$ is even. The quoted formulae simply perform the computation of the Euler characteristic from the combinatorial features of the representing graph.

As a particular case, when $m=3$ (with $d \ge 4$), we obtain the following formula, which holds for any $(d+1)$-colored graph $(\Gamma,\gamma)$ representing a singular $d$-manifold, with $d\ge 4$:
\begin{equation} \label{g_rst}
2g_{r,s,t}= g_{r,s} + g_{s,t} + g_{r,t} -p.
\end{equation}

Also relations \eqref{residue_d_odd} and \eqref{residue_d_even} directly follow from the computation of the Euler characteristic via combinatorial elements of the representing graph: in fact,   if  
$\mathcal B = \Delta_d - \{i\}$ and all $i$-colored vertices of $K(\Gamma)$ are not singular,  each connected component of $\Gamma_{\hat\imath}$ represents the $d$-sphere.   

\smallskip

The difference between relation \eqref{genereGamma} and relation \eqref{genereGamma_c}, by making use of relation \eqref{g_rst} applied to the subset $\mathcal B= \{i-1, i, i +1\}$ of $\Delta_n,$ yields \eqref{3-residuesversus2-residues}. 
On the other hand, the difference between relation \eqref{3-residuesversus2-residues} and the same relation applied to the graph $\Gamma_{\widehat {\varepsilon_j}}$ 
(for $j \notin \{i-1, i, i +1\}$) yields \eqref{numerospigoli}. 
  
\smallskip  

As a consequence, since $g_{\widehat{\varepsilon_{i}},\widehat{\varepsilon_{j}}} \ge  g_{\widehat{\varepsilon_{i}}} + g_{\widehat{\varepsilon_{j}}} -g$   trivially holds, we have: 
\begin{equation} 
\rho_\varepsilon(\Gamma) - \rho_{\varepsilon}(\Gamma_{\widehat {\varepsilon_i}}) - \rho_{\varepsilon}(\Gamma_{\widehat {\varepsilon_j}}) + \rho_{\varepsilon}(\Gamma_{\widehat {\varepsilon_j}, \widehat {\varepsilon_i}}) \ge 0  \ \ \ \ \forall i,j \text{\ non consecutive in \ } \varepsilon.
\end{equation}

Moreover, by applying formula  \eqref{numerospigoli} to the graph $\Gamma_{\widehat {\varepsilon_i}}$, we obtain: 
\begin{equation} \label{numerotriangoli}
g_{\widehat{\varepsilon_{i}},\widehat{\varepsilon_{j}}\widehat{\varepsilon_{k}}} =  (g_{\widehat{\varepsilon_{i}},\widehat{\varepsilon_{k}}} + g_{\widehat{\varepsilon_{i}},\widehat{\varepsilon_{j}}} -g_{\widehat{\varepsilon_{i}}} ) +  \rho_\varepsilon(\Gamma_{\widehat{\varepsilon_{i}}}) - \rho_{\varepsilon}(\Gamma_{\widehat {\varepsilon_i}, \widehat{\varepsilon_{k}}}) - \rho_{\varepsilon}(\Gamma_{\widehat {\varepsilon_i}, \widehat{\varepsilon_{j}}}) + \rho_{\varepsilon}(\Gamma_{\widehat {\varepsilon_i}, \widehat {\varepsilon_j}, \widehat{\varepsilon_{k}}})  \end{equation}
$\forall j,k \text{\ non consecutive in \ } (\varepsilon_{0},\ldots,\varepsilon_{i-1},\varepsilon_{i+1},\ldots,\varepsilon_{d}).$   
 \medskip

Now, if $j_1, j_2, \dots, j_{d-3} \notin \{i-1, i, i +1\},$ the difference between relation \eqref{3-residuesversus2-residues} and the same relation applied to the graph $\Gamma_{\widehat {\varepsilon_{j_1}}, \widehat {\varepsilon_{j_2}}, \dots \widehat {\varepsilon_{j_{d-3}}}}$ yields:
\begin{eqnarray} \label{3-residuesversus2-residues-bis}
&&g_{\widehat{\varepsilon_{i}},\widehat {\varepsilon_{j_1}}, \widehat {\varepsilon_{j_2}}, \dots \widehat {\varepsilon_{j_{n-3}}}} =  
(g_{\widehat{\varepsilon_{i}}} + g_{\widehat {\varepsilon_{j_1}}, \widehat {\varepsilon_{j_2}}, \dots \widehat {\varepsilon_{j_{d-3}}}} -g ) +\nonumber\\  
&&+ [\rho_\varepsilon(\Gamma) - \rho_{\varepsilon}(\Gamma_{\widehat {\varepsilon_i}}) - (\rho_{\varepsilon}(\Gamma_{\widehat {\varepsilon_{j_1}}, \widehat{\varepsilon_{j_2}}, \dots \widehat{\varepsilon_{j_{d-3}}}}) -
\rho_{\varepsilon}(\Gamma_{\widehat {\varepsilon_i}, \widehat{\varepsilon_{j_1}}, \widehat{\varepsilon_{j_2}}, \dots \widehat{\varepsilon_{j_{d-3}}}}))].
\end{eqnarray}

Note that $\{\widehat{\varepsilon_{i}},\widehat {\varepsilon_{j_1}}, \widehat {\varepsilon_{j_2}}, \dots \widehat {\varepsilon_{j_{d-3}}}\} = \{\varepsilon_{i-1}, \varepsilon_{i+1}, \varepsilon_{r}\}$ with $r \in \Delta_d - \{i-1, i, i +1\}$. 
Hence, the previous relation may be written as:
\begin{eqnarray*} g_{\varepsilon_{i-1}, \varepsilon_{i+1}, \varepsilon_{r}} =&&  (g_{\widehat{\varepsilon_{i}}} + g_{\varepsilon_{i-1}, \varepsilon_{i}, \varepsilon_{i+1}, \varepsilon_{r}} -g ) +  \rho_\varepsilon(\Gamma) - \rho_{\varepsilon}(\Gamma_{\widehat {\varepsilon_c}}) +\\
&&- \rho_{\varepsilon}(\Gamma_{\varepsilon_{c-1}, \varepsilon_{c}, \varepsilon_{c+1}, \varepsilon_{r}})
+\rho_{\varepsilon}(\Gamma_{\varepsilon_{i-1}, \varepsilon_{i+1}, \varepsilon_{r}}).\end{eqnarray*}

Since $\Gamma_{i-1,i+1,r}$ represents $\mathbb S^2$, $\rho_{\varepsilon}(\Gamma_{i-1,i+1,r})=0$ holds; hence, we may further simplify the relation as:
\begin{equation} \label{3-residuesversus2-residues-quater}
g_{\varepsilon_{i-1}, \varepsilon_{i+1}, \varepsilon_{r}} =  (g_{\widehat{\varepsilon_{i}}} + g_{\varepsilon_{i-1}, \varepsilon_{i}, \varepsilon_{i+1}, \varepsilon_{r}} -g ) +  \rho_\varepsilon(\Gamma) - \rho_{\varepsilon}(\Gamma_{\widehat {\varepsilon_i}}) - \rho_{\varepsilon}(\Gamma_{\varepsilon_{i-1}, \varepsilon_{i}, \varepsilon_{i+1}, \varepsilon_{r}}).
\end{equation}

Finally, by comparing relation \eqref{3-residuesversus2-residues} and relation \eqref{3-residuesversus2-residues-quater}, we obtain relation \eqref{2-residuesversuspairsof3-residues}. \qed

\bigskip

\begin{rem} {\em Note that relation \eqref{3-residuesversus2-residues} also yields:
\begin{equation}\label{chiavacci-pareschi}
\rho_\varepsilon(\Gamma) - \rho_{\varepsilon}(\Gamma_{\widehat {\varepsilon_i}}) = (g - g_{\widehat{\varepsilon_i}}) + (g_{\varepsilon_{i-1},\varepsilon_{i+1}} - g_{\varepsilon_{i-1}, \varepsilon_{i}, \varepsilon_{i+1}}) \geq 0\ \ \ \ \forall i \in \Delta_d.
\end{equation}
In \cite{Chiavacci-Pareschi}  the inequality $ \rho_{\varepsilon}(\Gamma_{\widehat {\varepsilon_i}}) \le \rho_\varepsilon(\Gamma) $ was already proved to hold for any $(d+1)$-colored graph.}
\end{rem}

\bigskip

\begin{prop} \label{equality-subgraphs}
Let $(\Gamma,\gamma)$ be a bipartite (resp. non-bipartite) $(d+1)$-colored graph representing a singular $d$-manifold. 
\begin{itemize}
\item If a color $i \in \Delta_d$ exists such that $\rho_\e (\Gamma) = \rho_{\varepsilon}(\Gamma_{\widehat{\varepsilon_i}})$, \
then 
$ \rho_{\varepsilon}(\Gamma_{\mathcal B}) = \rho_{\varepsilon}(\Gamma_{\mathcal B-\{i\}})$ \ 
for each subset $\mathcal B \subset \Delta_d$ with $\{i-1, i, i +1\} \subset \mathcal B.$ 

\item In particular, if $\mathcal B=\{i-1, i, i +1,r\}$, \ $\rho_\e (\Gamma) = \rho_{\varepsilon}(\Gamma_{\widehat{\varepsilon_i}})$ implies $\rho_{\varepsilon}(\Gamma_{\mathcal B})=0.$    
\end{itemize}
\end{prop}

\dimo
As a consequence of relations \eqref{numerospigoli} and \eqref{3-residuesversus2-residues-bis}, we have: 
\begin{eqnarray*} 0 &\le& \rho_{\varepsilon}(\Gamma_{\widehat {\varepsilon_{j_1}}, \widehat{\varepsilon_{j_2}}, \dots \widehat{\varepsilon_{j_{n-3}}}}) - 
\rho_{\varepsilon}(\Gamma_{\widehat {\varepsilon_i}, \widehat{\varepsilon_{j_1}}, \widehat{\varepsilon_{j_2}}, \dots \widehat{\varepsilon_{j_{n-3}}}}) \le \\  
&\le&\rho_{\varepsilon}(\Gamma_{\widehat {\varepsilon_{j_1}}, \widehat{\varepsilon_{j_2}}, \dots \widehat{\varepsilon_{j_{n-2}}}}) - 
\rho_{\varepsilon}(\Gamma_{\widehat {\varepsilon_i}, \widehat{\varepsilon_{j_1}}, \widehat{\varepsilon_{j_2}}, \dots \widehat{\varepsilon_{j_{n-2}}}}) \le \dots\\
 \dots  &\le&  \rho_{\varepsilon}(\Gamma_{\widehat {\varepsilon_{j_1}}, \widehat{\varepsilon_{j_2}}}) - \rho_{\varepsilon}(\Gamma_{\widehat {\varepsilon_i}, 
 \widehat{\varepsilon_{j_1}}, \widehat{\varepsilon_{j_2}}}) \le \rho_{\varepsilon}(\Gamma_{\widehat {\varepsilon_{j_1}}}) - \rho_{\varepsilon}(\Gamma_{\widehat {\varepsilon_{j_1}}, 
 \widehat {\varepsilon_i}}) 
\end{eqnarray*}

for each $j_1, j_2, \dots, j_{d-3} \notin \{i-1, i, i +1\}.$

The first  statement now easily follows. As regards the second one, it is sufficient to note that,   in case  $\mathcal B=\{i-1, i, i +1,r\}$, $\Gamma_{\mathcal B-\{i\}}$ represents a 2-dimensional sphere, and hence its regular genus is zero. \qed
 
\bigskip
\section{General properties of generalized regular genus} 
\label{sec_general-properties}

Within crystallization theory, two standard methods are known, in order to obtain a presentation of the fundamental group of a closed manifold directly from a graph representing it. The following extensions to compact manifolds and singular manifolds hold: 

\begin{prop}  \label{fundamental_group} 
Let $(\Gamma,\gamma)$ be a $(d+1)$-colored graph representing the singular $d$-manifold $N$ and the associated compact $d$-manifold $\check N$. 

\par \noindent   
\begin{itemize}
\item For each $i,j \in \Delta_d$,  let $X_{ij}$  (resp. $R_{ij}$) be a set in bijection with the connected components of $\Gamma_{\hat\imath \hat\jmath}$ (resp. with the $\{i,j\}$-colored cycles of $\Gamma$), and let $\bar R_{ij}$ be a subset  of $X_{ij}$ corresponding to the a maximal tree of the subcomplex $K_{ij}$ of $K(\Gamma)$ (consisting only of vertices labelled $i$ and $j$, and edges connecting them). 
Then: 
\begin{itemize}
\item[(a)]
if $i,j\in \Delta_d$ are not singular in $\Gamma$, \ \ \ $$  \pi_1(\check N) \ = \ < X_{ij} \ / \ R_{ij} \cup \bar R_{ij}>;$$ 
\item[(a')]
if no color in $\Delta_d - \{i,j\}$ is singular in $\Gamma$,  \ \ \ $$  \pi_1(N) \ = \ < X_{ij} \ / \ R_{ij} \cup \bar R_{ij} >.$$
\end{itemize}

\par \noindent
\item For each $i \in \Delta_d$, let $X_{i}$   (resp. $R_{i}$) be a set in bijection with the $i$-colored edges of $\Gamma$ (resp. with the $\{i,j\}$-colored cycles of $\Gamma$, for any $j \in \Delta_d - \{i\}$) and let $\bar R_{i}$  be a subset  of $X_{i}$ corresponding to a minimal set of $i$-colored edges of $\Gamma$ connecting $\Gamma_{\hat\imath}.$ Then: 
\begin{itemize}
\item[(b)]
if $i \in \Delta_d$ is not singular  in $\Gamma$, 
$$  \pi_1(\check N) \ = \ < X_{i} \ / \ R_{i} \cup \bar R_{i} >;$$   
\item[(b')]
if no color in $\Delta_d - \{i\}$ is singular in $\Gamma$, 
$$  \pi_1(N) \ = \ < X_{i} \ / \ R_{i} \cup \bar R_{i} >.$$   
\end{itemize}
\end{itemize}
\end{prop}

\dimo
It is a direct consequence of some general results concerning the fundamental groups of pseudocomplexes associated to colored graphs: see \cite{Chiavacci}.  \qed

\bigskip
The following statement yields an interesting inequality between the generalized regular genus and the rank of the fundamental group, for any compact manifold with connected boundary. The analogous inequality for closed manifolds is well known: see   \cite[Proposition B]{Chiavacci-Pareschi}.

\begin{prop}  \label{genus vs rank}
Let $M$ be a compact $d$-manifold with empty or connected boundary.  
Then: 
$$  \bar {\mathcal G}(M) \ge rk(\pi_1(M))$$
\end{prop}

\dimo
Let $(\Gamma, \gamma)$ be a $(d+1)$-colored graph realizing the generalized regular genus of $M$, 
with respect to the permutation $\varepsilon$ of $\Delta_d$, i.e.   $ {\rho_{\varepsilon}}(\Gamma)=  \bar {\mathcal G}(M)$. 
Let $i$ and $j$ be two not singular colors that are not consecutive in the permutation $\varepsilon$: they certainly exist since $M$ has empty or connected boundary and so $\Gamma$  has at most one singular color. 
It is now sufficient to consider the presentation of the fundamental group of $M$ given by  Proposition  \ref{fundamental_group}(a), with respect to colors $i$ and $j$ and to recall that, in virtue of formulae \eqref{numerospigoli} and \eqref{chiavacci-pareschi},
$$ \# (X_{ij} -  \bar R_{ij}) \le g_{\widehat{\varepsilon_{i}},\widehat{\varepsilon_{j}}} -  (g_{\widehat{\varepsilon_{i}}} + g_{\widehat{\varepsilon_{j}}} -1 ) \le     \rho_\varepsilon(\Gamma).$$   \qed

\bigskip

Let us now recall  that another graph-based representation theory for compact (PL) manifolds exists, making use of colored graphs which fail to be regular. More precisely, any compact $d$-manifold can be represented by a pair $(\Lambda, \lambda)$, where $\lambda$ is still an edge-coloration on $E(\Lambda)$ by means of $\Delta_d$, but $\Lambda$ may miss some (or even all) $d$-colored edges: such a pair  
is said to be a {\it $(d+1)$-colored graph with boundary, regular with respect to color $d$}, and vertices missing the $d$-colored edge are called {\it boundary vertices} (see \cite{Ferri-Gagliardi-Grasselli}).

An easy combinatorial procedure, called {\it capping off}, enables to connect this representation to the one - involving only regular colored graphs - considered in Section \ref{prelim}.

\begin{prop}[\cite{Ferri-Gagliardi Yokohama 1985}] 
\label{cappingoff} Let $(\Lambda, \lambda)$ be a $(d+1)$-colored graph with boundary, regular with respect to color $d$, representing the compact $d$-manifold $M$.  
Chosen a color $c \in \Delta_{d-1}$, let $(\Gamma,\gamma)$ be the regular $(d+1)$-colored graph obtained from $\Lambda$ by {\rm capping off with respect to color $c$}, i.e. by joining  two boundary vertices by a $d$-colored edge, whenever they belong to the same $\{c,d\}$-colored path in $\Lambda.$ Then,  $(\Gamma,\gamma)$ represents the singular $d$-manifold $\widehat M$, and hence $M$, too.
\end{prop}

By means of (non-regular) edge-colored graphs with boundary, together with a suitable extension of Proposition \ref{reg_emb}, 
Gagliardi introduced within crystallization theory a ``classical" notion of {\it regular genus} for compact $d$-manifolds, too (see \cite{Gagliardi 1981} and \cite{Gagliardi_boundary}).
The following result establishes a comparison between regular genus and generalized regular genus (as defined in Section \ref{prelim}: see Definitions \ref{Gurau-degree} and \ref{def_gen_degree}) for any compact $d$-manifold.

\begin{prop}  \label{regular genus vs generalized regular genus}
Let $M$ be a compact $d$-manifold, with $d \ge 3$, and let $\mathcal G(M)$ denote the regular genus of $M$.
Then:   
$$  \overline{\mathcal G}(M) \ \le \ \mathcal G(M),$$
Moreover:
\begin{itemize}
\item[(a)] $\overline{\mathcal G}(M) \  = \ \mathcal G(M)$   if $M$ is a closed $d$-manifold;
\item[(b)] $\overline{\mathcal G}(M) \  = \ 0 \ \iff  \  M \ = \ \mathbb S^d  \ \iff  \ \mathcal G(M)\ = \ 0$;   
\item[(c)] $\overline{\mathcal G}(M) \  = \ \mathcal G(M)$   if $M$ is a compact $3$-manifold with connected boundary;   
\item[(d)] there exist compact $3$-manifolds (with disconnected boundary) such that \ $\overline{\mathcal G}(M) \  < \ \mathcal G(M).$  
\end{itemize}

\end{prop}
\dimo
The general inequality is a consequence of the ``capping off" procedure, recalled in Proposition \ref{cappingoff}. In fact, let us assume the regular genus of $M$ to be realized by the (not regular) graph with boundary $\Lambda$ with respect to the cyclic permutation $\varepsilon = (\e_0, \e_1, \dots, \e_{d-1}, \e_d=d)$ of $\Delta_d$. Then, it is not difficult to prove that, if $c\in \{\e_0, \e_{d-1}\}$ is chosen, 
and $\Gamma$ is the (regular) $(d+1)$-colored graph obtained from $\Lambda$ by capping off with respect to color $c$, the generalized regular genus of $\Gamma$ with respect to $\e$ equals the regular genus of $\Lambda$  with respect to the same permutation: $\rho_{\e} (\Gamma)=\rho_{\e} (\Lambda) = \mathcal G(M).$ 

Equality (a) is trivial by definition (as already pointed out in Remark \ref{rem_regular genus}). 

Regarding statement (b), first note that, obviously, $\overline{\mathcal G}(\mathbb S^d)=\mathcal G(\mathbb S^d)=0$; moreover, the main theorem of \cite{Ferri-Gagliardi Proc AMS 1982} ensures that, if $\Gamma$ represents a closed $d$-manifold $M$,  $\rho(\Gamma)=0$ implies $M$ to be a PL $d$-sphere. In order to complete the proof of both co-implications, let us consider a (regular) $(d+1)$-colored graph $\Gamma$ such that there exists a cyclic permutation $\varepsilon$ of $\Delta_d$ with $\rho_\varepsilon(\Gamma) = 0$; we want to prove that $|K(\Gamma)|$ is a closed $d$-manifold. 
If $d=2$ then $|K(\Gamma)|\cong\mathbb S^2$, since $\rho_\varepsilon(\Gamma)$ trivially coincides with the genus of the surface $|K(\Gamma)|.$
Suppose now our claim to be true in each dimension $<d$; given $i\in\mathbb Z_{d+1}$, let $\Xi$ be a connected component of $\Gamma_{\widehat{\varepsilon_i}}$, which is a $d$-colored graph.
Since $\rho_{\varepsilon}(\Xi)\leq\rho_\varepsilon(\Gamma)$ (see inequality \eqref{chiavacci-pareschi})  then, by induction, $\Xi$ represents a PL $(d-1)$-sphere and,
therefore, $|K(\Gamma)|$ is a closed PL $d$-manifold, and statement (b) is proved. 

The inequality $\overline{\mathcal G}(M) \ \ge \ \mathcal G(M)$ for $3$-dimensional manifolds with connected boundary, yielding relation (c), is proved in \cite{Cristofori-Mulazzani}.  

The same paper also presents examples of the strict inequality (d): if $F$ is a closed surface of genus $g$, then \ $\overline{\mathcal G}(F \times I)=g \ < \ \mathcal G(F \times I)=2g.$
\qed

\begin{rem} \label{generalized vs regular genus}
{\em In Section \ref{section.4dim} we will prove that the equality between the two invariants does not hold for $4$-manifolds with boundary, even if the boundary is assumed to be connected: see Corollary \ref{cor.no-equality (n=4)}(b).}
\end{rem}

As regards the invariant regular genus, a well-known relation (i.e. $\mathcal G(M)  \ge \, \mathcal G(\partial M)$) compares the regular genus of any compact manifold with the regular genus of its boundary; in the case of connected boundary, the following extensions hold, concerning both the generalized regular genus and the G-degree: 
  
\begin{prop}  \label{genus vs boundary genus}
Let $M$ be a compact  $d$-manifold with (non-empty) connected boundary. Then: 
$$  \bar {\mathcal G}(M) \, \ge \, \mathcal G(\partial M) \quad \quad \text{and} \quad \quad  \mathcal D_G(M) \, \ge  \, d \cdot \mathcal D_G(\partial M). $$
\end{prop}

\dimo
The first inequality is an easy consequence of \eqref{chiavacci-pareschi}, applied to a regular graph $\Gamma$ representing $M$, so that $\bar {\mathcal G}(M) =\rho_\varepsilon(\Gamma)$ ($\varepsilon$ being a cyclic permutation of $\Delta_d$) and having color $i$ as its (only) singular color. 

The second inequality may be obtained in a similar way, by making use of the relation $\omega_G (\Gamma) \, \ge \, d \cdot  \omega_G (\Gamma_{\hat\imath}),$ proved in \cite[Lemma 4.6]{Gurau-Ryan} for each $(d+1)$-colored graph and for each color $i \in \Delta_d.$ 
\qed

A $d$-dimensional extension of the construction described in \cite[Proposition 5(i)]{Cristofori-Mulazzani} (resp. in \cite[Proposition 5(ii)]{Cristofori-Mulazzani}), performed in \cite[Section 7]{Grasselli-Mulazzani} in a general setting  including
graphs representing singular $d$-manifolds,  allows to easily obtain graphs representing connected sums (resp. boundary connected sums) of  compact PL $d$-manifolds directly from the graphs representing the summands. 

We briefly recall that, 
if $(\Gamma_1,\gamma_1)$ and $(\Gamma_2,\gamma_2)$ are two disjoint $(d+1)$-colored graphs and $v_i \in V_i$ for each $i \in \{1,2\},$  the \emph{graph connected sum} of $\Gamma_1$, $\Gamma_2$  with respect to vertices $v_1, v_2$ (denoted by $\Gamma_1\#_{v_1v_2}\Gamma_2$)  
is the graph obtained from $\Gamma_1$ and $\Gamma_2$ by deleting $v_1$ and $v_2$ and welding the ``dangling'' edges of the same color. It is not difficult to check that, if all $d$-residues containing  $v_1$ and $v_2$ 
are ordinary (resp. if both $v_1$ and $v_2$ belong to exactly one singular $d$-residue, $\Xi_1$ and $\Xi_2$, say), then  $\Gamma_1\#_{v_1v_2}\Gamma_2$ represents the (internal) connected sum between $|K(\Gamma_1)|$ and 
$|K(\Gamma_2)|$  (resp. represents the boundary  connected sum between $|K(\Gamma_1)|$ and $|K(\Gamma_2)|$,  performed on the boundary components corresponding to $\Xi_1$ and $\Xi_2$).

\begin{prop}  \label{connected sums}
Let $M_1$ and $M_2$ be compact $d$-manifolds.  Then: 
$$\bar{ \mathcal G}(M_1  \# M_2) \, \le \, \bar{ \mathcal G}(M_1) \, + \, \bar{ \mathcal G}(M_2) \quad \text{and} \quad  
\bar{ \mathcal G}(M_1  \,^{\partial} \# \, M_2) \, \le \, \bar{ \mathcal G}(M_1) \, + \, \bar{ \mathcal G}(M_2);$$  
$${ \mathcal D_G}(M_1  \# M_2) \, \le \, { \mathcal D_G}(M_1) \, + \, { \mathcal D_G}(M_2) \quad \text{and}$$ 
$${ \mathcal D_G}(M_1  \,^{\partial} \# \, M_2) \, \le \, { \mathcal D_G}(M_1) \, + \, { \mathcal D_G}(M_2).$$  \end{prop}

\dimo
It is an easy consequence of the above described constructions: see the quoted papers, together with  \cite[Proposition 10]{ Casali-Cristofori-Dartois-Grasselli}.
\qed

\section{Representing handlebodies and products $ \times \, I$} \label{sec_hadlebodies and xI}

\begin{prop}\label{prop.handlebodies}
For any $d \ge 4$, a bipartite (resp. non-bipartite) $(d+1)$-colored graph exists, with order $2d$ and regular genus one with respect to any permutation of $\Delta_d$, representing the genus one $d$-dimensional handlebody  $\mathbb Y^d_1$ (resp.  $\tilde{  \mathbb Y}^d_1$).
Hence, for each $d \ge 4:$
\begin{equation} \label{Y^d_1}
\bar{ \mathcal G}(\mathbb Y^d_1)=   \bar{ \mathcal G}(\tilde{\mathbb Y}^d_1) = 1 \quad \quad \text{and} \quad \quad \mathcal D_G (\mathbb Y^d_1)=  \mathcal D_G (\tilde{\mathbb Y}^d_1) = \frac{d!}2.
\end{equation} 
Moreover, for each $d \ge 4$ and for each $m \ge 1:$   
\begin{equation}  \label{Y^d_m}
\bar{ \mathcal G}(\mathbb Y^d_m)  =   \bar{ \mathcal G}(\tilde{\mathbb Y}^d_m)  =m \quad \quad \text{and} \quad \quad \mathcal D_G (\mathbb Y^d_m)=  \mathcal D_G (\tilde{\mathbb Y}^d_m) = m  \cdot \frac{d!}2. \end{equation} 
\end{prop}

\begin{figure}
\centering 
\includegraphics[scale=.5]{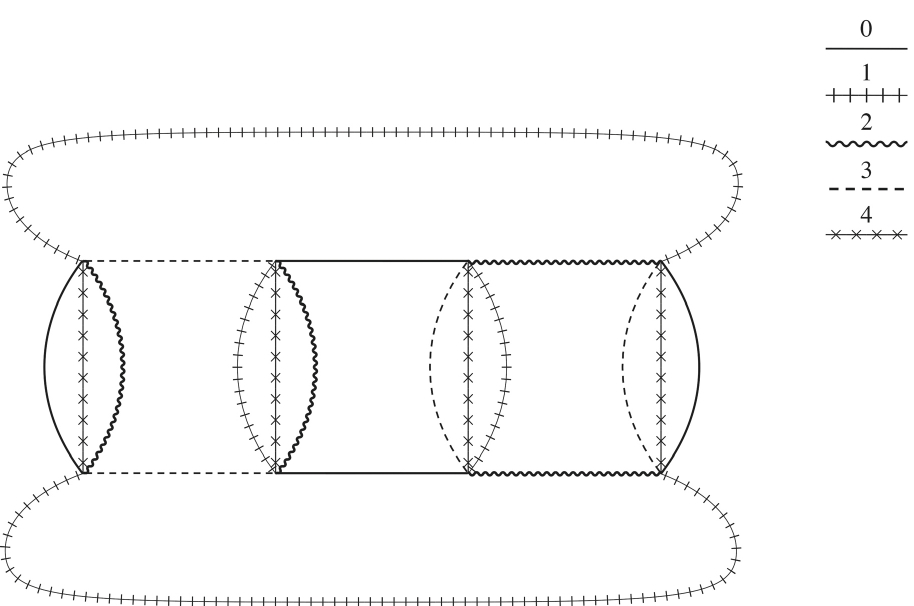}
\caption{Regular 5-colored graph representing ${\mathbb Y}^4_1.$}
\label{fig1:handlebody}
\end{figure}

\begin{figure}
\centering
\includegraphics[scale=.6]{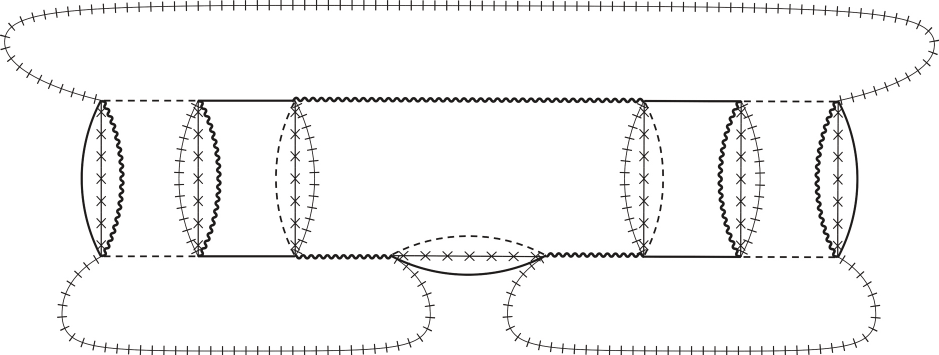}
\caption{Regular 5-colored graph representing the boundary connected sum ${\mathbb Y}^4_1 \,^{\partial} \# {\mathbb Y}^4_1 = {\mathbb Y}^4_2.$}
 \label{Figure_boundary-connected-sum}
\end{figure}

\dimo
For any $d \ge 3$, an order $2(d+1)$  $(d+1)$-colored graph with boundary $(H,h)$ (resp. $(H^{\prime},h^{\prime})$) is well-known, which represents the genus one $d$-dimensional handlebody  $\mathbb Y^d_1$ (resp.  $\tilde{\mathbb Y}^d_1$) 
: see \cite{Gagliardi_boundary}.  
By applying to $(H,h)$ (resp. $(H^{\prime},h^{\prime})$) the ``capping off" procedure described  in Proposition \ref{cappingoff}, a (regular) order $2(d+1)$  $(d+1)$-colored graph representing $\mathbb Y^d_1$ (resp.  $\tilde{ \mathbb Y}^d_1$) is obtained. It is easy to check that it admits a (proper) 2-dipole, whose elimination yields a (minimal) order $2d$ regular $(d+1)$-colored graph $(\hat{H},\hat{h})$ (resp. $(\hat{H^{\prime}},\hat{h^{\prime}})$) representing $\mathbb Y^d_1$ (resp.  $\tilde{  \mathbb Y}^d_1$): see Figure \ref{fig1:handlebody} for the orientable $4$-dimensional case.  A direct computation gives  $\rho_{\varepsilon}(\hat{H})=1$ (resp.   $\rho_{\varepsilon}(\hat{H}^{\prime})=1$) for each permutation $\varepsilon$ of $\Delta_d.$ 
Hence, the classification of compact PL $d$-manifolds with generalized  regular genus zero (and with G-degree zero) easily allows to prove \eqref{Y^d_1}: see Proposition \ref{regular genus vs generalized regular genus} (b).     
Now, it is not difficult to check that, for each $m \ge 1$, the graph connected sum construction hinted to in the previous Section, with suitable choices of the vertices, enables to obtain a bipartite (resp. non-bipartite) $(d+1)$-colored graph representing  the genus $m$ $d$-dimensional handlebody  $\mathbb Y^d_m = \, ^\partial \#_m  \mathbb Y^d_1$ (resp.  $\tilde{  \mathbb Y}^d_m = \,^\partial \#_m  \tilde{\mathbb Y}^d_1$): its order is $2md -2(m-1)$  and its regular genus is $m$ with respect to any permutation of $\Delta_d$. See Figure \ref{Figure_boundary-connected-sum} for an example, in case $d=4$ and $m=2$. 
As a consequence of this construction, together with the inequalities $ \bar{ \mathcal G}(M) \ge \mathcal G (\partial M)$ and  ${ \mathcal D_G}(M) \ge d \cdot \mathcal D_G (\partial M)$ (see Proposition  \ref{genus vs boundary genus}),  both equalities of  \eqref{Y^d_m} easily follow. 
\qed

\medskip

\medskip

\begin{figure}
\centering
\includegraphics[scale=.6]{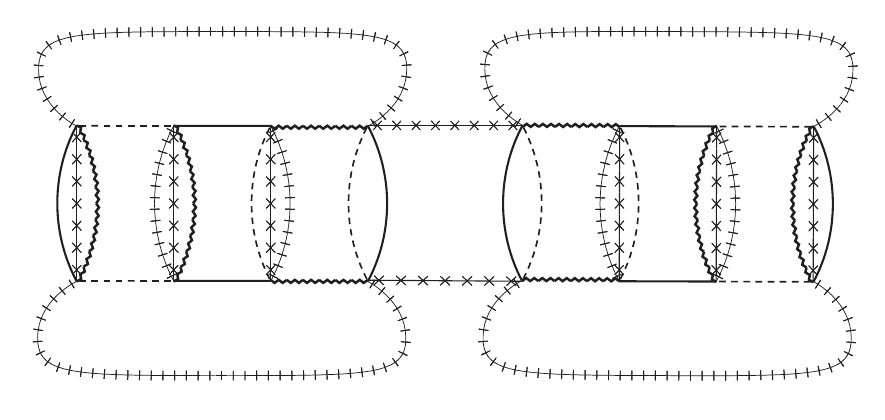}
\caption{Regular 5-colored graph representing the connected sum ${\mathbb Y}^4_1 \# {\mathbb Y}^4_1.$}
 \label{Figure_connected-sum}
\end{figure}

\begin{rem} {\em Note that, by a suitable application of the procedure of graph connected sum,  it is easy to obtain also a (bipartite or non-bipartite) $(d+1)$-colored graph representing the connected sum of $m$ ($m \ge 2$)  (orientable or non-orientable) $d$-dimensional  handlebodies. 
If  $\tilde{\tilde{\mathbb Y}}^d_{r}$ denotes either the orientable or non-orientable genus $r$ $d$-dimensional handlebody, then the graph representing $\tilde{\tilde{\mathbb Y}}^d_{r_1}  \#  \cdots \# \tilde{\tilde{\mathbb Y}}^d_{r_m}$ has order $2d(r_1 + \cdots + r_m)$  (since the procedure has to be preceded by the insertion of $m$ $d$-dipoles, in order to obtain $m$ ordinary $d$-residues)  and its regular genus is $r_1 + \cdots + r_m$ with respect to any permutation of $\Delta_d$.  See Figure \ref{Figure_connected-sum}  for an example, in case $d=4$, $m=2$ and $r_1=r_2=1$.  \\ 
\\ The following inequalities directly follow by construction:\footnote{Note that, in virtue of  relation  $ \mathcal D_G(M) \ \ge \ \frac{d!}2 \cdot \overline{\mathcal G}(M)$ (see Section \ref{prelim}), if the equality concerning generalized regular genus holds, then the one concerning G-degree holds, too.} }
\begin{equation}   \label{mY^d_1}
\bar{ \mathcal G}(\tilde{\tilde{\mathbb Y}}^d_{r_1}  \#  \cdots \# \tilde{\tilde{\mathbb Y}}^d_{r_m})   \le r_1 + \cdots + r_m  \quad \quad 
\mathcal D_G (\tilde{\tilde{\mathbb Y}}^d_{r_1}  \#  \cdots \# \tilde{\tilde{\mathbb Y}}^d_{r_m}) \le (r_1 + \cdots + r_m) \cdot \frac{d!}2. \end{equation} 
\end{rem}

\begin{prop}\label{MxI}
For any $d \ge 3$, given a bipartite (resp. non-bipartite) $(d+1)$-colored graph $(\Lambda, \lambda)$ representing a closed orientable (resp. non-orientable) $d$-manifold $M$, then a bipartite (resp. non-bipartite) $(d+2)$-colored graph $(\tilde \Lambda, \tilde \lambda)$ representing the singular $(d+1)$-manifold ${M \times I}$ exists, with $\#V(\tilde \Lambda)= \#V(\Lambda)$ and $\rho(\tilde \Lambda)= \rho(\Lambda)$. 
Hence, for each $d \ge 3:$
$$ \bar{ \mathcal G}(M \times I) =  {\mathcal G}(M)  \ \ \text{for any closed $d$-manifold} \ M.$$
Moreover:
$$\omega_G(\tilde \Lambda)  \ = \ \frac{d!}2  \Big[\sum_{i \in \{1, \dots d\}} (p - g_{\varepsilon_0 \varepsilon_i}) - (p-1) + \frac 2 {(d-1)!} \omega_G(\Lambda)\Big],$$
$\varepsilon = (\varepsilon_0 ,\varepsilon_1, \dots ,\varepsilon_d)$ being the cyclic permutation of $\Delta_d$ so that $\rho(\Lambda)= \rho_{\varepsilon}(\Lambda)$.  
\end{prop}

\dimo
Let $\varepsilon = (\varepsilon_0 ,\varepsilon_1, \dots ,\varepsilon_d)$ be the cyclic permutation of $\Delta_d$ so that $\rho(\Lambda)= \rho_{\varepsilon}(\Lambda)$. 
If $(\tilde \Lambda, \tilde \lambda)$ is obtained from $\Lambda$ by adding a $(d+1)$-colored edge between any pair of $\varepsilon_0$-adjacent vertices, then it is easy to check that $\tilde \Lambda$ represents $M \times I$ and $\rho_{\varepsilon^{\prime}}(\tilde \Lambda) = \rho(\tilde \Lambda),$ where $\varepsilon^{\prime} = (\varepsilon_0 ,\varepsilon_1, \dots,\varepsilon_d, d+1):$ see  \cite{Ferri-Gagliardi Yokohama 1985}, and Figure \ref{Fig_L(2,1)xI} for an example of the construction, with $M=L(2,1)$.
This fact proves the first part of the statement and, as a consequence, the inequality $ \bar{ \mathcal G}(M \times I) \le  {\mathcal G}(M)$ for any closed $d$-manifold $M.$  
On the other hand, since  $M \times I$ has two boundary components PL-homeomorphic to $M$, any $(d+2)$-colored graph representing  $M \times I$ as a singular $(d+1)$-manifold must have a $(d+1)$-residue representing $M$, and hence must have regular genus greater or equal to ${\mathcal G}(M)$. 
As regards the computation of the G-degree, a direct application of formula  $ \omega_G(\Gamma) \ = \ \frac {(d-1)!} 2 \Big(d + \frac d 2 (d-1) p - \sum_{r,s\in \Delta_d} g_{rs} \Big) $ (proved in \cite{Casali-Cristofori-Dartois-Grasselli} for any $(d+1)$-colored graph)  to $(\tilde \Lambda, \tilde \lambda)$ yields (in virtue of the combinatorial structure of $\tilde \Lambda$): 
\begin{eqnarray*}\omega_G(\tilde \Lambda)  \ =&& \ \frac{d!}2  \Big[(d+1) + \frac {(d+1) dp} 2 -  \sum_{r,s\in \Delta_{d+1}} g_{\varepsilon_r \varepsilon_s} \Big] \ = 
\\ = && \ \frac{d!}2  \Big[(d+1) + \frac {(d+1) dp} 2 -  \sum_{r,s\in \Delta_{d} }g_{\varepsilon_r \varepsilon_s} -  \sum_{i\in \{1, \dots d\}}  g_{\varepsilon_0 \varepsilon_i} - p \Big].
\end{eqnarray*}

Moreover, by making use of the similar computation for  $ \omega_G(\Lambda)$, yielding $\sum_{r,s\in \Delta_{d}} g_{\varepsilon_r \varepsilon_s} = d +  \frac d 2 (d-1) p -  \frac 2 {(d-1)!} \omega_G(\Lambda)$, we obtain: 
\begin{eqnarray*} \omega_G(\tilde \Lambda)  \ =&& \ \frac{d!}2  \Big[(d+1) + \frac {(d+1) dp} 2 -  d - \frac d 2 (d-1) p -  \sum_{i\in \{1, \dots d\}}  g_{\varepsilon_0 \varepsilon_i} - p\ + \\ &&+  \frac 2 {(d-1)!} \omega_G(\Lambda) \Big],
\end{eqnarray*}
which the last formula of the statement follows from. 
\qed

\begin{figure}
\centering
\includegraphics[scale=.5]{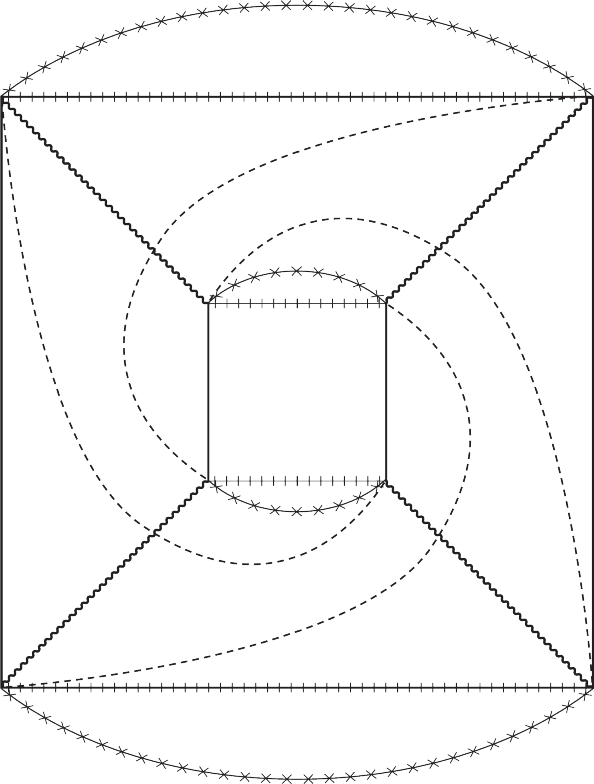}
\caption{Regular 5-colored graph representing the product $L(2,1) \times I.$}
 \label{Fig_L(2,1)xI}
\end{figure}

\section{Representing  $\mathbb D^2$-bundles over $\mathbb S^2$}
\label{section_bundles}

In  \cite{Casali UMI},  a $5$-colored graph with boundary $(\tilde \Gamma_c, \tilde \gamma_c)$  (resp. $(\tilde \Gamma_0, \tilde \gamma_0)$)  representing  the $\mathbb D^2$-bundle over $\mathbb S^2$,  $\xi_c$, with Euler class $c$  (resp.  the trivial  $\mathbb D^2$-bundle over $\mathbb S^2$,   $\mathbb S^2 \times \mathbb D^2$) is produced,  $\forall c \in \mathbb  Z^+ - \{1\}$; all these graphs have regular genus equal to three. This allows to prove - by means also of some theoretical results about the ``gap" between regular genus and the rank of the fundamental group of a PL 4-manifold with boundary - that $ {\mathcal G}(\xi_c)= {\mathcal G}(\mathbb S^2 \times \mathbb D^2)=3.$

The regular $5$-colored graphs obtained from the above graphs by means of the ``capping off" procedure described  in Proposition \ref{cappingoff} 
(which represent the singular 4-manifolds $\hat \xi_c$ and $\widehat{\mathbb S^2 \times \mathbb D^2}$, and hence the compact 4-manifolds $\xi_c$ and $\mathbb S^2 \times \mathbb D^2,$ too) have regular genus three by construction, have the same order as the starting graphs with boundary (i.e. $4c + 6$ for $\xi_c$, $\forall c \in \mathbb  Z^+ - \{1\},$ and $14$ for $\mathbb S^2 \times \mathbb D^2$),  but admit a (proper) 2-dipole involving colors non-consecutive in the permutation $\varepsilon$ realizing the minimum generalized regular genus, together with two 2-dipoles involving colors consecutive in the permutation $\varepsilon$. 

\begin{figure}
\centering
\includegraphics[scale=.55]{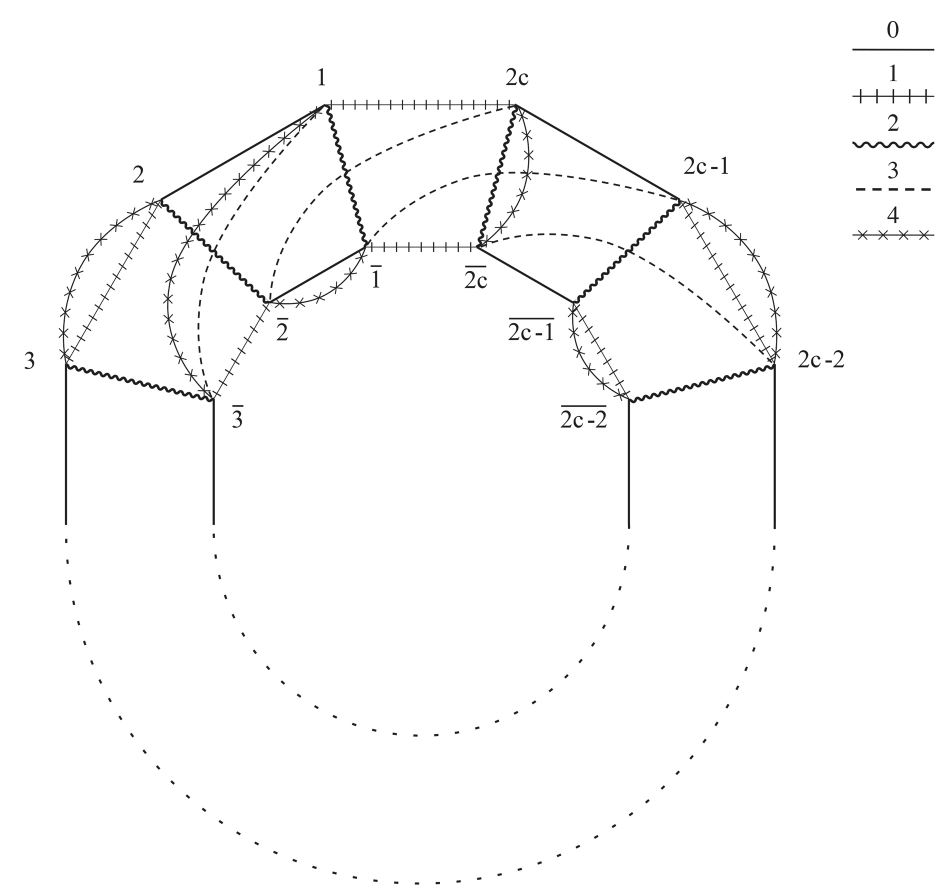}
\caption{The regular 5-colored graph $\Lambda_c$ representing the compact 4-manifold $\xi_c.$}
 \label{fig.csi_c}
\end{figure}

\begin{figure}
\centering
\includegraphics[scale=.5]{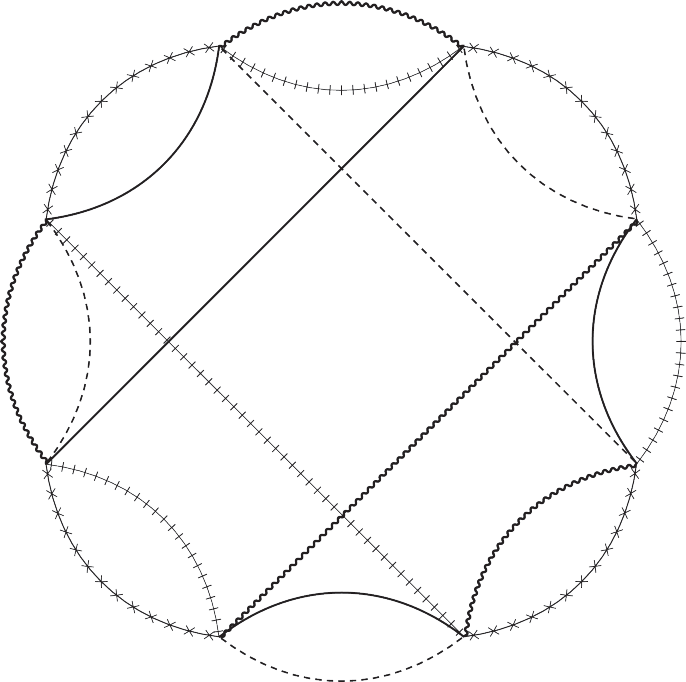}
\caption{The regular 5-colored graph $\Lambda_0$ representing the compact 4-manifold $\mathbb S^2 \times \mathbb D^2.$}
 \label{fig.S2xD2}
\end{figure}

Now, it is easy to check, via Proposition \ref{reg_emb} and  Definition \ref{def_dipole}, that the elimination of a 2-dipole involving colors non-consecutive (resp. consecutive) in the permutation $\varepsilon$ decreases by one (resp. does not affect)  the regular genus with respect to  $\varepsilon$. 
Hence, the elimination of the three 2-dipoles yields a regular 5-colored graph $(\Lambda_c, \lambda_c)$, $\forall c \in \mathbb Z^+ - \{1\}$ (resp. $(\Lambda_0, \lambda_0)$)  representing $\xi_c$ (resp. $\mathbb S^2 \times \mathbb D^2$) with the same order $4c$ (resp. $8$) as the standard crystallization of $L(c,1)$ (resp. of  \, $\mathbb S^1 \times \mathbb S^2$): see Figure \ref{fig.csi_c}  (resp. Figure \ref{fig.S2xD2}).     

As a consequence, we  have:  
$$ \bar{ \mathcal G}(\xi_c) \le 2   \ \ \forall c \in \mathbb  Z^+ - \{1\} \ \ \ \ \text{and} \ \ \ \ \ \bar {\mathcal G}(\mathbb S^2 \times \mathbb D^2) \le 2.$$

Actually, in the following Corollary \ref{D2-fibrati}, we will prove that all compact 4-manifolds of this infinite class turn out to have generalized regular genus equal to two. 

\medskip 
As far as the G-degree is concerned, we recall that  \cite[Proposition 5]{Casali-Grasselli 2017} proves, for each $5$-colored graph representing a compact PL $4$-manifold:   
\begin{equation}\label{multiplo6}\omega_G(\Gamma)= 6 \Big(\rho_\varepsilon(\Gamma) + \rho_{\varepsilon^{\prime}}(\Gamma)\Big),\end{equation}
where $\varepsilon=(\varepsilon_0, \varepsilon_1, \varepsilon_2, \varepsilon_3, \varepsilon_4)$ is an arbitrary permutation  of $\Delta_4$ and $ \varepsilon^{\prime}$ is the  ``associated" permutation, i.e.  $ \varepsilon^{\prime} = (\varepsilon_0, \varepsilon_2, \varepsilon_4, \varepsilon_1, \varepsilon_3)$.

A direct computation allows to check that,  if $\varepsilon^{\prime}$  denotes the permutation associated to the one realizing $\rho(\Lambda_c)=2$ (resp.  $\rho(\Lambda_0)=2$), then  the regular genus with respect to $\varepsilon^{\prime}$  is $2c-2$ $\forall c \in \mathbb  Z^+ - \{1\}$ (resp. is also $2$)
Then, formula \eqref{multiplo6} yields $\omega_G(\Lambda_c)= 6[2+(2c-2)]=12c$ \ \  $\forall c \in \mathbb  Z^+ - \{1\}$ \  (resp. $\omega_G(\Lambda_0)= 6 (2+2) =24$). 

Hence: 
$$ \mathcal D_G(\xi_c) \le 12c  \ \ \forall c \in \mathbb  Z^+ - \{1\} \ \ \ \ \text{and} \ \ \ \ \ \mathcal D_G(\mathbb S^2 \times \mathbb D^2) \le 24.$$

\section{Classifying results in dimension 4}  
\label{section.4dim}

\subsection{Classifying with respect to generalized regular genus}  
\label{subsection.4dim-genus}

In the 4-dimensional setting, formula \eqref{g_rst} enables to prove the following useful results.

\begin{prop} \label{prop-Euler-characteristic n=4}
Let $(\Gamma,\gamma)$ be a connected $5$-colored graph representing a singular 4-manifold $N^4$ (and the associated compact  $4$-manifold $\check N^4$). 
For each cyclic permutation $\varepsilon$ of $\Delta_4:$ 
\begin{equation} \label{Euler-caracteristic n=4}
 \chi (N^4)  = 2 - 2 \rho_{\varepsilon}(\Gamma) + \sum_{i \in \Delta_4} \rho_{\varepsilon}(\Gamma_{\widehat{ \varepsilon_i}}).
\end{equation}

Moreover, if  $\rho_\e$ and  $\rho_{\widehat{\e_{i}}}$ respectively denote   $\rho_\e(\Gamma)$ and $\rho_{\varepsilon}(\Gamma_{\widehat{ \varepsilon_i}}):$ 
\begin{equation} \label{numerospigoli n=4}
g_{\widehat{\e_{i-1}},\widehat{\e_{i+1}}} = g_{\e_i,\e_{i+2}, \e_{i+3}} \ =  \   (g_{\widehat {\e_{i-1}}}  +  g_{\widehat {\e_{i+1}}}  -1) +     \rho_\e - \rho_{\widehat{\e_{i-1}}} - \rho_{\widehat{\e_{i+1}}};  
 \end{equation} 
 \begin{equation} \label{numerotriangoli n=4}
g_{\widehat{\e_i},\widehat{\e_j},\widehat{\e_k}} =  (g_{\widehat{\e_i},\widehat{\e_j}} + g_{\widehat{\e_i},\widehat{\e_k}} -g_{\widehat{\e_i}})  +  \rho_{\widehat{\e_i}}  \ \ \ \ \ \ \forall j,k \text{\ non consecutive in \ } \Delta_4-\{i\};  \end{equation}
\begin{equation} \label{sum_subgenus n=4}
\sum_{i \in \Delta_4} g_{\widehat{\e_{i-1}},\widehat{\e_{i+1}}} = \Big( 2 \sum _{i \in \Delta_4}  g_{\hat {\e_i}}   -5 \Big) + 5 \rho_\e - 2  \sum_{i \in \Delta_4}  \rho_{\widehat{\e_{i}}};
 \end{equation}  
\begin{equation} \label{inequality_two-subgenus}
\rho_{\widehat{\e_{i-1}}} + \rho_{\widehat{\e_{i+1}}} \le \rho_\e;  
 \end{equation}
 \begin{equation} \label{inequality_subgenus}
\sum_{i \in \Delta_4}  \rho_{\widehat{\e_i}} \le \Big\lfloor \frac 5 2 \, \rho_\e \Big \rfloor;
\end{equation}
\begin{equation} \label{equality_subgenus}
\sum_{i \in \Delta_4}  \rho_{\widehat{\e_i}} =  \frac 5 2 \, \rho_\e  \ \ \ \ \Longleftrightarrow \ \ \ \  \rho_{\widehat{\e_{i-1}}} + \rho_{\widehat{\e_{i+1}}} = \rho_\e \ \ \ \ \forall i \in \Delta_4.
\end{equation}
 \end{prop}

\dimo
Relations \eqref{numerospigoli n=4} and  \eqref{numerotriangoli n=4} are nothing but  relations \eqref{numerospigoli} and \eqref{numerotriangoli}, in case $\Gamma$ is assumed to be a connected graph representing a singular 4-manifold. 

Summing up relations  \eqref{numerospigoli n=4}, for each $i \in \Delta_4$, yields relation  \eqref{sum_subgenus n=4}. 

Since $g_{\widehat{\e_{i-1}},\widehat{\e_{i+1}}} \ge g_{\widehat {\e_{i-1}}}  +  g_{\widehat {\e_{i+1}}}  -1$ trivially holds, relations \eqref{inequality_two-subgenus} and  \eqref{inequality_subgenus}  follow from relations \eqref{numerospigoli n=4} and  \eqref{sum_subgenus n=4} respectively. 

Finally, the co-implication \eqref{equality_subgenus} is a direct consequence.
\qed

\begin{cor} \label{genere-subgeneri (sfera n=4)}
Let $(\Gamma,\gamma)$ be a connected $5$-colored graph representing $\mathbb S^4.$ Then, for each cyclic permutation $\varepsilon$ of $\Delta_4$:  
\begin{equation} \label{genere-subgeneri(sfera n=4)}
\rho_{\varepsilon}(\Gamma) = \frac 1 2 \sum_{i \in \Delta_4} \rho_{\varepsilon}(\Gamma_{\widehat{ \varepsilon_i}}).
\end{equation} \end{cor}
\ \qed

\medskip

\begin{prop} \label{null-subgenus}
Let $(\Gamma,\gamma)$ be a bipartite (resp. non-bipartite) $5$-colored graph representing a compact PL $4$-manifold $M^4$ with empty or connected boundary  
and let $\e$ be a cyclic permutation of $\Delta_4.$
If there exists  $i \in \Delta_4$ so that $\rho_{\varepsilon}(\Gamma_{\widehat{ \varepsilon_i}})=0,$
\ then  $M^4 \cong \#_{\alpha} (\mathbb S^1 \times \mathbb S^3) \# \mathbb Y^4_{\beta}$  (resp. 
$M^4 \cong \#_{\alpha} (\mathbb S^1 \tilde \times \mathbb S^3)\# {\mathbb Y}^4_{\beta}$ or $M^4 \cong \#_{\alpha} (\mathbb S^1 \tilde \times \mathbb S^3)\# \tilde{\mathbb Y}^4_{\beta}$)
with  $\alpha, \beta \ge 0$ and $\alpha + \beta \le \rho_\e (\Gamma)$ and 
$\beta \le \rho_{\varepsilon}(\Gamma_{\widehat{ \varepsilon_c}}),$ $c\in \Delta_4$ being the singular color of $\Gamma$ (if any\ \footnote{If all colors are not singular, then $\beta=0$ as in case (a).}).

In particular: 
\begin{itemize}
\item[(a)] if $M^4$ is a closed $4$-manifold and there exists  $i \in \Delta_4$ so that $\rho_{\varepsilon}(\Gamma_{\widehat{ \varepsilon_i}})=0,$  
then $M^4 \cong \#_{\alpha} (\mathbb S^1 \times \mathbb S^3)$  (resp. $M^4 \cong \#_{\alpha} (\mathbb S^1 \tilde \times \mathbb S^3)$), with  $\alpha \le \rho_\e (\Gamma);$ 
\item[(b)] 
if $M^4$ has (non-empty) connected boundary,  and $\ \rho_\e (\Gamma) = \rho_{\varepsilon}(\Gamma_{\widehat{ \varepsilon_c}})$, \ $c \in \Delta_4$ being the only singular color of $\Gamma$, 
\ then  $M^4 \cong \mathbb Y^4_{m}$  (resp. $M^4 \cong \tilde{\mathbb Y}^4_{m}$), with $m \le \rho_\e (\Gamma)$.
\end{itemize}
\end{prop}

\dimo
Let $\varepsilon$ be the cyclic permutation of $\Delta_4$ such that $ \rho (\Gamma)= \rho_{\e}(\Gamma)$; without loss of generality we may assume $\varepsilon=(0,1,2,3,4)$ and, by Remark \ref{rem_dipoli},  
$g_{\hat\imath}=1$ for each $i \in \Delta_4.$    
Moreover, since the cyclic permutation $\e$ is defined up to inverse, we may further assume that $i \in \{c+1, c+2\}$, \ $c\in \Delta_4$ being the color of the possible singular vertex.

Now, relation  \eqref{numerotriangoli n=4} yields: 
 $$ g_{\widehat{{c-1}}, \widehat{{c+1}}, \widehat{{c+2}}} =
g_{\widehat{{c-1}}, \widehat{\imath}}  + g_{\widehat{{c+1}}, \widehat{{c+2}}} - g_{\hat\imath}=1 + \rho_{\varepsilon}(\Gamma_{\widehat{ \varepsilon_i}})  = g_{\widehat{{c-1}}, \widehat{\imath}}  + g_{\widehat{{c+1}}, \widehat{{c+2}}}  -1.$$
Hence, $K(c-1,c+1,c+2)$ collapses onto $K(c-1,j)$, where $\{j \}=  \{c+1, c+2\} - \{i\}$. Since $c-1,c+1,c+2$ are not singular colors, in the bipartite (resp. non-bipartite)  case $N(c-1,c+1,c+2) = \mathbb Y^4_m$           
(resp. $N(c-1,c+1,c+2) = \tilde{\mathbb Y}^4_m$) follows, with $m= g_{\widehat{{c-1}}, \widehat{\jmath}} -1$.

On the other hand,  $K(c,c-2)$ consists of  $g_ {\widehat{{c}}, \widehat{{c-2}}}$ edges; since $c-2$ is not a singular color, $N(c,c-2)$ is obtained by the cone over $lkd (v_c)$ by adding $g_ {\widehat{{c}}, \widehat{{c-2}}} -1 $ 1-handles. 

Since $\widehat M^4=  N(c-1,c+1,c+2) \cup_\phi N(c,c-2),$  $\phi$ being a boundary homeomorphism, in the bipartite (resp. non-bipartite)  case  $\partial N(c,c-2) = \partial N(c-1,c+1,c+2) =  \partial \mathbb Y^4_m  
=  \#_{m} (\mathbb S^1 \times \mathbb S^2)$ (resp. $\partial N(c,c-2) = \partial N(c-1,c+1,c+2) =  \partial \tilde{\mathbb Y}^4_m         
=  \#_{m} (\mathbb S^1 \tilde \times \mathbb S^2)$) follows. 

\noindent In the orientable case,  we have $lkd (v_c) \#   [\#_{m^{\prime}} (\mathbb S^1 \times \mathbb S^2)] = \#_{m} (\mathbb S^1 \times \mathbb S^2)$, where  $m^{\prime}= g_{\widehat{{c}}, \widehat{{c-2}}}-1$, i.e.  (in virtue of the uniqueness of the sum decomposition in dimension 3) $lkd (v_c) =   \#_{m-m^{\prime}} (\mathbb S^1 \times \mathbb S^2).$ 
The statement now follows from  \cite[Lemma 1]{Casali-Malagoli}, since: 
\begin{eqnarray*}\widehat M^4 &=&  [\ {\mathbb Y}^4_{m^{\prime}} \ ^{\partial} \#  \  {\mathbb Y}^4_{m-m^{\prime}} ]  \cup_\phi  [{\mathbb Y}^4_{m^{\prime}}  \ ^{\partial}  \# \  [v_c * \#_{m-m^{\prime}} (\mathbb S^1 \times \mathbb S^2)] ]  = \\
&=&  [  {\mathbb Y}^4_{m^{\prime}} \cup_\phi   {\mathbb Y}^4_{m^{\prime}} ]  \#  [ {\mathbb Y}^4_{m-m^{\prime}}   \cup_\phi   [v_c * \#_{m-m^{\prime}} (\mathbb S^1 \times \mathbb S^2)] ] = \\
&=& \#_{m^{\prime}} (\mathbb S^1  \times \mathbb S^3) \# \hat  { {\mathbb Y}}^4_{m-m^{\prime}}. \end{eqnarray*}
                                                                                                               
The non-orientable case may be proved in full analogy, by distinguishing the case of orientable and non-orientable boundary.

\medskip

Both point (a) and point (b) of the statement are nothing but particular cases of the general statement; indeed, as regards point (b), the second part of Proposition  \ref{equality-subgraphs} yields, for $d=4$: $$  \rho_\e (\Gamma) = \rho_{\varepsilon}(\Gamma_{\widehat{ \varepsilon_c}})  
\ \ \Rightarrow \ \ \rho_{\varepsilon}(\Gamma_{\widehat {c+2}})  = \rho_{\varepsilon}(\Gamma_{\widehat {c-2}})  =0.$$
\qed

\begin{rem} {\em Note that point (a) of the above proposition could be independently proved simply by noting that relation \eqref{numerotriangoli n=4} yields $$g_{i,i+2} = g_{i, i+1, i+2} + g_{i-1,i, i+2} -1,$$ which ensures - via Lemma 5 of \cite{Casali - Forum Math 1992} - that $K(i-2,i-1,i+1)$ collapses to a graph, i.e.  $N(i-2,i-1,i+1)$ is a handlebody. 
Since also $N(i,i+2)$ is obviously a handlebody, statement (a) follows via a well-known theorem by Montesinos and Laudenbach-Poenaru (see \cite{Montesinos} and \cite{Laudenbach-Poenaru}). \\
Also point (b) could be independently proved by noting that, by formula \eqref{numerospigoli n=4}, $g_{\hat{c},\widehat{c+2}} =1$ follows, 
i.e. $N(c,c+2)$ is homeomorphic to the cone over $lkd(v_c).$ Moreover, relation \eqref{numerotriangoli n=4} yields 
$$ g_{\widehat{c-2},\widehat{c-1}, \widehat{c+1}} = g_{\widehat{c-2},\widehat{c-1}} +  g_{\widehat{c+1},\widehat{c-2}} -1,$$ which ensures - via Lemma 5 of \cite{Casali - Forum Math 1992} - that $K(c-2,c-1,c+1)$ collapses to a graph, i.e.  $N(c-2,c-1,c+1)$ is a handlebody; this proves statement (b).}
\end{rem}

\bigskip

We are now able to classify all compact 4-manifolds with generalized regular genus one.

\begin{prop} \label{gen-genus_one}
Let $M^4$ be a compact 4-manifold with \ $\bar{ \mathcal G}(M^4)=1.$ 
\ Then, 
$$ \text{either}  \ \ M^4 \in \{\mathbb S^1 \times \mathbb S^3, \   \mathbb S^1 \tilde \times \mathbb S^3\}  
\ \ \ \text{or} \ \ \ M^4 \in \{{\mathbb Y}^4_1, \tilde{  \mathbb Y}^4_1\}   
\ \ \ \text{or} \ \ \  M^4 \cong \bar M \times I, $$

where $\bar M$ is a genus one closed 3-manifold.  
\end{prop}

\dimo   
Three cases occur: 
\begin{itemize}
\item $M^4$ is a closed 4-manifold; 
\item  $M^4$ is a compact 4-manifold with (non-empty) connected boundary; 
\item $M^4$ is a  compact 4-manifold with disconnected boundary. 
\end{itemize}

In the first and second case, if  $\Gamma$ represents $M^4$ with $\rho_{\varepsilon}(\Gamma)= \rho(\Gamma) = \bar{ \mathcal G}(M^4) = 1,$ then $\Gamma$ may be assumed to satisfy $g_{\hat\imath}=1$ $\forall i \in \Delta_4$  (see Remark \ref{rem_dipoli}). 
Relation \eqref{inequality_subgenus} directly implies the existence of at least a color $i \in \Delta_4$ such that $\rho_{\varepsilon}(\Gamma_{\widehat{ \varepsilon_i}})=0.$ Hence,  Proposition \ref{null-subgenus} proves the statement: see in particular points (a) and (b).   

Let us take into account the third case, i.e. $K(\Gamma)$, with $|K(\Gamma)|= \widehat M^4$, contains more than one singular vertex. 
By using the notations of Proposition \ref{prop-Euler-characteristic n=4},  if $c$ is one of the colors of singular vertices, $\rho_{\hat{c}} \le \rho =1$ obviously implies that $g_{\hat{c}} =1$ may be assumed to hold, i.e. $K(\Gamma)$ contains only one $c$-colored singular vertex.  Let now $d$ be the color of another singular vertex: by relation  \eqref{inequality_two-subgenus}, $d \in \{c-1, c+1\}$ follows, while relation  \eqref{inequality_subgenus} implies $ \rho_{\hat{\imath}} =0$ $\forall i \in \Delta_4 -\{c, d\}$ (i.e. $K(\Gamma)$ contains exactly two singular vertices). 
Without loss of generality, let us assume $\rho_{\hat{c}} = \rho_{\widehat{c+1}} = \rho =1$  and   $ \rho_{\hat{\imath}} =0$ $\forall i \in \Delta_4 -\{c, c+1\}.$ 

Since $g_{\hat{c},\widehat{c+2}} =  1 + \rho - \rho_{\hat{c}} - \rho_{\widehat{c+2}} = 1+1-1-0=1,$  $K(c,c+2)$ consists of one only edge, i.e. $N(c,c+2)$ is homeomorphic to the cone over $lkd(v_c)$.

On the other hand, by relation \eqref{numerotriangoli n=4}, $\rho_{\widehat{c-1}}=0$ yields  \ $ g_{\widehat{c-2},\widehat{c-1}, \widehat{c+1}} = g_{\widehat{c-2},\widehat{c-1}} +  g_{\widehat{c+1},\widehat{c-2}} -1;$
hence - via Lemma 5 of \cite{Casali - Forum Math 1992} - $K(c-2,c-1,c+1)$ is proved to collapse onto $K(c-2,c+1)$, which consists of exactly one edge (since $g_{\widehat{c+1},\widehat{c-2}} = 1 + (\rho - \rho_{\widehat{c-2}}) - \rho_{\widehat{c+1}}= 1 + 1 - 0 -1 =1$). This proves that  $N(c-2,c-1,c+1)$ is homeomorphic to the cone over $lkd(v_{c+1})$, too. 

Thus $M^4 \cong \bar M \times I$, where $\bar M$ is a genus one closed 3-manifold (homeomorphic to both $lkd(v_c)$ and $lkd(v_{c+1})$), now easily follows. 
\qed

\begin{cor} \label{D2-fibrati}
Let  $\xi_c$ be the $\mathbb D^2$-bundle over $\mathbb S^2$ with Euler class $c$, $\forall c \in \mathbb  Z^+ - \{1\}$. Then, 
$$ \bar{ \mathcal G}(\xi_c)= \bar {\mathcal G}(\mathbb S^2 \times \mathbb D^2)=2.$$
Moreover, 
$$\bar{ \mathcal G}(\mathbb Y^4_1 \# \mathbb Y^4_1)  =  \bar{ \mathcal G}(\mathbb Y^4_1 \# \tilde{ \mathbb Y}^4_1)  =  \bar{ \mathcal G}(\tilde{ \mathbb Y}^4_1 \# \tilde{ \mathbb Y}^4_1)  = 2. $$
\end{cor}

\dimo
By Proposition \ref{gen-genus_one}, $ \bar{ \mathcal G}(\xi_c) \ge 2$, $\bar {\mathcal G}(\mathbb S^2 \times \mathbb D^2) \ge 2$,  $\bar{ \mathcal G}(\mathbb Y^4_1 \# \mathbb Y^4_1)  \ge 2 $, $ \bar{ \mathcal G}(\mathbb Y^4_1 \# \tilde{ \mathbb Y}^4_1)  \ge 2$  and  $\bar{ \mathcal G}(\tilde{ \mathbb Y}^4_1 \# \tilde{ \mathbb Y}^4_1) \ge 2$ trivially follow.

On the other hand, in Section \ref{section_bundles} (resp. in Section \ref{sec_hadlebodies and xI}), we have obtained 5-colored graphs with generalized regular genus two representing $\mathbb S^2 \times \mathbb D^2$ and $\xi_c$, $\forall c \in \mathbb  Z^+ - \{1\}$ (resp. representing  $\mathbb Y^4_1 \# \mathbb Y^4_1$,  $\mathbb Y^4_1 \# \tilde{ \mathbb Y}^4_1$ and  $\tilde{ \mathbb Y}^4_1 \# \tilde{ \mathbb Y}^4_1$): see Figures \ref{fig.S2xD2} and \ref{fig.csi_c} (resp. see Figure \ref{Figure_connected-sum}).
Hence, the statement is proved.
\qed

The results about non-finiteness-to-one of generalized regular genus (already pointed out in Section \ref{intro} and in Remark \ref{generalized vs regular genus}) now easily follow:    

\begin{cor}  \label{cor.no-equality (n=4)}   
\par \noindent
\begin{itemize}
\item[(a)] 
Generalized regular genus is not finite-to-one in dimension four. 
\item[(b)]
In dimension four,  the equality between regular genus and generalized regular genus of manifolds with boundary does not hold, even if the boundary is assumed to be connected. 
\end{itemize}
\vskip-0.5truecm
\ \qed
\end{cor}

Further results are obtained,  concerning compact PL $4$-manifolds with generalized regular genus two.     

\begin{prop} \label{gen-genus_two}
Let $M^4$ be a compact 4-manifold  with empty or connected boundary, 
with \ $\bar{ \mathcal G}(M^4)=2.$ 
\ Then: 
\begin{itemize}  
\item[$\bullet$]  either  $M^4 \in \{\#_2 (\mathbb S^1 \times \mathbb S^3),  \#_2(\mathbb S^1 \tilde \times \mathbb S^3),  \mathbb{CP}^2 \},$
\item[$\bullet$]  or  $M^4 \in \{{\mathbb Y}^4_2,  \tilde{\mathbb Y}^4_2, 
(\mathbb S^1 \times \mathbb S^3)\# {\mathbb Y}^4_1,  (\mathbb S^1 \tilde \times \mathbb S^3)\# {\mathbb Y}^4_1,  (\mathbb S^1 \tilde \times \mathbb S^3)\# \tilde{\mathbb Y}^4_1,
\mathbb S^2 \times \mathbb D^2,$ $\xi_2\}, $ 
\item[$\bullet$]  or   $M^4 \cong M^4(K,d),$ \ $ (K,d)$ being a framed knot such that $M^3(K,d) =  L(\alpha, \beta)$ with $\alpha \ge 3.$
\end{itemize}  
\end{prop}

\dimo
Let  $(\Gamma,\gamma)$ be  a $5$-colored graph representing $M^4$, with  $\rho (\Gamma)= \bar{ \mathcal G}(M^4)=2$.  
Without loss of generality, for sake of simplicity we may assume that the cyclic permutation $\e$ of $\Delta_4$ such that $ \rho (\Gamma)= \rho_{\e}(\Gamma)$ is $\varepsilon=(0,1,2,3,4)$ and that $g_{\hat\imath}=1$ holds for each $i \in \Delta_4$ (see Remark \ref{rem_dipoli}).  

In virtue of the inequality  \eqref{inequality_two-subgenus},  $\rho (\Gamma)= 2$ implies that the only possible cases are: 
\begin{itemize}
\item there exists $i \in \Delta_4$ so that  \ $ \rho_{\hat\imath} =0$; 
\item  $\rho_{\hat{\imath}} = 1$ \ $\forall i \in \Delta_4.$ 
\end{itemize}

In the first case, by Proposition \ref{null-subgenus}   \   either $M^4 \cong \#_{\alpha} (\mathbb S^1 \times \mathbb S^3) \# \mathbb Y^4_{\beta}$  or  
 $M^4 \cong \#_{\alpha} (\mathbb S^1 \tilde \times \mathbb S^3)\# {\mathbb Y}^4_\beta$ or $M^4 \cong \#_{\alpha} (\mathbb S^1 \tilde \times \mathbb S^3)\# \tilde{\mathbb Y}^4_\beta$ 
hold, with $ \alpha, \beta \ge 0$ and $\alpha + \beta = 2$. 
Hence, if $M^4$ is a closed 4-manifold, $M^4 \in \{\#_2 (\mathbb S^1 \times \mathbb S^3),  \#_2(\mathbb S^1 \tilde \times \mathbb S^3)\}$ follows,  while $M^4 \in \{{\mathbb Y}^4_2, \tilde{\mathbb Y}^4_2,$ 
$(\mathbb S^1 \times \mathbb S^3) \# \mathbb Y^4_1, (\mathbb S^1 \tilde \times \mathbb S^3)\# {\mathbb Y}^4_1, (\mathbb S^1 \tilde \times \mathbb S^3)\# \tilde{\mathbb Y}^4_1 \}$  
 follows if $M^4$ has connected boundary.

\medskip
Let us now take into account the second case. If $c$ is the color of the only (possible) singular vertex $v_c$ of $K(\Gamma)$ or any color if $M^4$ is closed, 
relations  \eqref{numerospigoli n=4} and \eqref{numerotriangoli n=4}  yield: 
$$  g_{\widehat{i-1},\widehat{i+1}} =  1+ \rho - \rho_{\widehat{i-1}} - \rho_{\widehat{i+1}} = 1 \ \ \forall i \in \Delta_4 
\ \ \ \ 
\text{(in particular: \ } g_{\hat{c},\widehat{c-2}} =  1 )$$ 
and  $$g_{\widehat{{c-1}}, \widehat{{c+1}}, \widehat{{c+2}}} =
g_{\widehat{{c-1}}, \widehat{{c+1}}}  + g_{\widehat{{c+1}}, \widehat{{c+2}}} + \rho_{\widehat{c+1}} 
- 1 = g_{\widehat{{c-1}}, \widehat{{c+1}}}  + g_{\widehat{{c+1}}, \widehat{{c+2}}}= 1 +  g_{\widehat{{c+1}}, \widehat{{c+2}}}.$$
Hence, $M^4$ is simply-connected (in virtue of Proposition  \ref{fundamental_group}) and $N(c,c-2)$ is  the cone over $lkd (v_c),$  while $K(c-1,c+1,c+2)$ collapses - by arguments already used in Lemma 2 of \cite{Casali UMI} - to two 2-simplices with common boundary, i.e. $N(c-1,c+1,c+2)$ is obtained from a 0-handle $H^{(0)}= \mathbb D^4$ by addition of one 2-handle $H^{(2)}$, according to a framed knot $(K,d).$    

Now, if $v_c$ is not singular 
(that is, if $M^4$ is a closed $4$-manifold),   $N(c,c-2)$ is a 4-dimensional disk, and hence - by a well-known theorem in \cite{Gordon-Luecke} - $(K,d)$ turns out to be the trivial knot with framing $1$. So, $M^4 \cong \mathbb{CP}^2$ directly follows. 
 
On the other hand, if $v_c$ is a singular vertex, $M^4 \cong M^4(K,d)$ holds, \ $(K,d)$ being a framed knot such that $M^3(K,d)$ has genus one (equal to $\rho_{\hat{c}}$). 

If $M^3(K,d)\cong \mathbb S^1 \times \mathbb S^2$, a classic result of Dehn surgery ensures $(K,d)$ to be the 0-framed trivial knot (see \cite{Gabai}), i.e. $M^4 \cong \mathbb S^2 \times \mathbb D^2.$
Further, if $M^3(K,d)\cong L(2,1)$, another, more recent, result of Dehn surgery ensures $(K,d)$ to be the 2-framed trivial knot (see \cite{Kronheimer}),  i.e. $M^4 \cong \xi_2.$

Hence, the only remaining cases concern simply-connected $4$-manifolds $M^4(K,d)$ having lens spaces $L(\alpha, \beta)$, with $\alpha \ge 3$, as boundary.  
\qed

\medskip

We are now able to prove the theorem, already stated in Section \ref{intro}, that summarizes the obtained classification results for compact $4$-manifolds according to their generalized regular genus.

\medskip
\dimo (Theorem \ref{Thm.gen-genus_one&two})

Statement (a) is nothing but the case $d=4$ of Proposition \ref{regular genus vs generalized regular genus}(b).   

Statement (b) is a direct consequence of Propositions \ref{prop.handlebodies}, \ref{MxI} and \ref{gen-genus_one}, together with the well-known existence of $5$-colored graphs of
regular genus one representing the two $\mathbb S^3$-bundles over $\mathbb S^1.$

With regard to statement (c), the result comes directly from Proposition \ref{gen-genus_two}, since for each $c\in\mathbb Z$, $\xi_c$, the $\mathbb D^2$-bundle over $\mathbb S^2$ with Euler class    
$c$, is exactly $M^4(K_0,c)$, $(K_0,c)$ being the $c$-framed trivial knot.
\qed

\begin{prop} \label{gen-genus_two_onesingularcolor}
Let  $(\Gamma,\gamma)$ be  a $5$-colored graph with exactly one singular color and  with  $\rho (\Gamma)= 2$.  
\ Then, either $K(\Gamma)$ has exactly one singular vertex  \! (and therefore $\Gamma$ represents one of the compact $4$-manifolds detected in Proposition \ref{gen-genus_two}), or 
$M^4 \in \{{\mathbb Y}^4_1 \#  {\mathbb Y}^4_1, \mathbb Y^4_1 \#  \tilde{ \mathbb Y}^4_1, \tilde{ \mathbb Y}^4_1 \#  \tilde{ \mathbb Y}^4_1\}.$ 
\end{prop}

\dimo 
Let $c$ be the singular color of $\Gamma$ and let $\e$ be the cyclic permutation of $\Delta_4$ such that $ \rho (\Gamma)= \rho_{\e}(\Gamma)=2$; further, let us assume - without loss of generality - $\varepsilon=(0,1,2,3,4)$ and $g_{\hat\imath}=1$ for each $i \in \Delta_4- \{c\}$ (see Remark \ref{rem_dipoli}). 
It is easy to check that, if $K(\Gamma)$ has more than one singular vertex, $\Gamma$ may be assumed to have exactly two $\hat c$-residues, both with regular genus one with respect to the permutation induced by $\e.$ Hence, $\rho_{\hat c}=2$.  Arguments similar to those used in the proof of Proposition  \ref{null-subgenus}(b)  ensure that $K(c,c+2)$ consists of two edges, with a common end-point (i.e. the $(c+2)$-labelled vertex) and with the other end-points consisting in the two singular $c$-labelled vertices of $K(\Gamma),$ $v^{\prime}_c$ and $v^{\prime \prime}_c$, say.  This easily implies that  $N(c,c+2)$ is homeomorphic to the boundary connected sum of a 4-disk, the cone  $v^{\prime}_c * lkd(v^{\prime}_c)$ and the cone $v^{\prime \prime}_c * lkd(v^{\prime \prime}_c)$; hence, the boundary of $N(c,c+2)$ is \ $lkd(v^{\prime}_c) \# lkd(v^{\prime \prime}_c)$.

On the other hand, formula \eqref{numerospigoli n=4} yields  $\rho_{\widehat{c-2}}=0$, and hence $ g_{\widehat{c-2},\widehat{c-1}, \widehat{c+1}} = g_{\widehat{c-2},\widehat{c-1}} +  g_{\widehat{c+1},\widehat{c-2}} -1,$ which implies that $K(c-2,c-1,c+1)$ collapses to a graph, i.e.  $N(c-2,c-1,c+1)$ is a handlebody of genus $m=g_{\widehat {c-1} \widehat {c+1}} -1 = 2- \rho_{\widehat {c-1}} - \rho_{\widehat {c+1}} \le 2$, whose boundary is a connected sum of $m\le 2$ (orientable or non-orientable) sphere bundles over $\mathbb S^1$.

Since the boundaries of $N(c,c+2)$  and $N(c-2,c-1,c+1)$ have to be identified, then $m=2$ and both $lkd(v^{\prime}_c)$ and $lkd(v^{\prime \prime}_c)$ must be homeomorphic to an orientable or non-orientable sphere bundle over $\mathbb S^1$. 

This easily proves that $M^4 \in \{{\mathbb Y}^4_1 \#  {\mathbb Y}^4_1, \mathbb Y^4_1 \#  \tilde{ \mathbb Y}^4_1, \tilde{ \mathbb Y}^4_1 \#  \tilde{ \mathbb Y}^4_1\},$ according to the orientability of the connected components of the boundary (i.e. of the singularities $lkd(v^{\prime}_c)$ and $lkd(v^{\prime \prime}_c)$).
\qed

\begin{thm}  \label{genus vs boundary genus (n=4)}
Let $M^4$ be a compact 4-manifold with (non-empty) connected boundary.  
Then: 
$$\bar{ \mathcal G}(M^4) = {\mathcal G}(\partial M^4) =m \ge 1   \quad  \Longleftrightarrow \quad   M^4 \in \{{\mathbb Y}^4_{m}, \tilde{  \mathbb Y}^4_m\} $$
\end{thm}

\dimo
Let $(\Gamma, \gamma)$ be a 5-colored graph realizing the generalized regular genus of $M^4$, with respect to the permutation $\varepsilon$ of $\Delta_4$: $ {\rho}(\Gamma) ={\rho_{\varepsilon}}(\Gamma) =m$. If $v_c$ is the only singular vertex of $|K(\Gamma)|=\widehat M^4$,  the subgraph $\Gamma_{\hat c}$ represents $\partial M^4$, and hence ${\rho_{\varepsilon}}(\Gamma_{\hat c}) \ge  {\mathcal G}(\partial M^4) =m.$ Since ${\rho_{\varepsilon}}(\Gamma_{\hat c}) \le {\rho_{\varepsilon}}(\Gamma)$ holds for any colored graph and for any color $c$, the equality  $ {\rho}(\Gamma) = {\rho_{\varepsilon}}(\Gamma_{\hat c})$ follows. 
     
$M^4 \in \{{\mathbb Y}^4_{m}, \tilde{  \mathbb Y}^4_m\} $        is now a direct consequence of   Proposition  \ref{null-subgenus}(b).  
Proposition \ref{prop.handlebodies} yields the converse implication of the statement.  
\qed

\medskip

\begin{thm}  \label{genus vs rank (n=4)} 
Let $M^4$ be a compact 4-manifold with empty or connected boundary. Then: 
\begin{itemize} 
\item[a.] \  $\bar{ \mathcal G}(M^4) \ge rk(\pi_1(M^4))$;
\item[b.] \   $\bar{ \mathcal G}(M^4)=rk(\pi_1(M^4))= \rho$ \ if and only if 
\begin{itemize}  
\item[$\bullet$] if $M^4$ is orientable,
\begin{itemize}  
\item[-] either  \  $M^4  \cong  \#_{\rho}(\mathbb S^1 \times \mathbb S^3)$
\item[-] or \  $M^4  \cong  \#_{\alpha} (\mathbb S^1  \times \mathbb S^3) \#  {\mathbb Y}^4_{\beta}$ with $\alpha + \beta = \rho$,
\end{itemize}
 \item[$\bullet$] if $M^4$ is non-orientable,
\begin{itemize}  
\item[-]  either \ $M^4   \cong  \#_{\rho}({\mathbb S^1 \, \widetilde{\times} \, \mathbb S^3})$
\item[-] or \    $M^4 \cong  \#_{\alpha} (\mathbb S^1  \widetilde{\times} \mathbb S^3) \# {\mathbb Y}^4_{\beta}$ or   
$M^4 \cong  \#_{\alpha} (\mathbb S^1  \widetilde{\times} \mathbb S^3) \#  \widetilde{\mathbb Y}^4_{\beta}$ again  with $\alpha + \beta  = \rho$.
\end{itemize}
\end{itemize}
\item[c.] \   $ \bar{ \mathcal G}(M^4) \ne rk(\pi_1(M^4)) \ \ \ \Longrightarrow \ \ \  \bar{ \mathcal G}(M^4)\,- \,rk(\pi_1(M^4)) \ \ge \ 2.$
\end{itemize}
\end{thm}

\dimo
The first statement is nothing but a particular case of Proposition  \ref{genus vs rank}. 
The proof of that result, in the 4-dimensional setting, yields (via formula \eqref{numerospigoli n=4}): 
$$ \# (X_{ij} - \bar R_{ij}) \le g_{\hat{\imath},\hat{\jmath}} - (g_{\hat {\imath}}  +  g_{\hat {\jmath}}  -1)   =    \rho - \rho_{\hat{\imath}} - \rho_{\hat{\jmath}},$$
for any pair $i,j$ of colors non-consecutive  in $ \varepsilon$ and such that both $\Gamma_{\hat\imath}$ and $\Gamma_{\hat\jmath}$ represent spheres.    
Hence, if $(\Gamma, \gamma)$ is a $5$-colored graph realizing the generalized regular genus of $M^4$, with respect to the permutation $\varepsilon$ of $\Delta_4$ (i.e. $ {\rho_{\varepsilon}}(\Gamma) = \rho =  \bar{ \mathcal G}(M^4)$), 
$\bar{ \mathcal G}(M^4)=rk(\pi_1(M^4))$ trivially implies  $ \rho_{\hat{\imath}} = \rho_{\hat{\jmath}}=0$, while $\bar{ \mathcal G}(M^4) - rk(\pi_1(M^4)) =1$ trivially implies  $ \rho_{\hat{\imath}} + \rho_{\hat{\jmath}} \le 1$. 

The second and third statements now easily follow from Proposition \ref{null-subgenus}, since all the represented compact $4$-manifolds actually satisfy the equality between the generalized regular genus and the rank of the fundamental group.
\qed

\subsection{Classifying with respect to G-degree}  
\label{subsection.4dim-degree}

In order to face the classification  
problem for compact $4$-manifolds with respect to G-degree, we need a further definition\footnote{Note that Definition \ref{gem-complexity} naturally extends - via graphs representing singular $d$-manifolds - an important PL invariant originally defined for closed manifolds. 
Many significant classification results have been obtained within crystallization theory with respect to it:  see, for example, \cite{BCrG} and \cite{Casali-Cristofori JKTR 2008} for the dimension $3$, \cite{Casali-Cristofori ElecJComb 2015} and \cite{Casali-Cristofori-Gagliardi Complutense 2015} for the dimension $4$. A classification 
according to gem-complexity for compact orientable $3$-manifolds with toric boundary is obtained in \cite{Cristofori-Fomynikh-Mulazzani-Tarkaev}.} and some preliminary results. 

\begin{defn} \label{gem-complexity}
{\em For each compact $d$-manifold $M$, its \emph{(generalized) gem-complexity} is the non-negative integer $k(M)= p - 1$, 
where $2p$ is the minimum order of a (regular) $(d+1)$-colored graph representing $M$.}
\end{defn}

\begin{lemma}\label{sommarho_i}
Let $(\Gamma,\gamma)$ be an order $2p$ 5-colored graph representing a compact $4$-manifold $M^4$. Then, for each cyclic permutation $\e$ of $\Delta_4$:
$$p = \frac 16 \omega_G(\Gamma) + 2\rho_\e(\Gamma) - \sum_{i \in \Delta_4} \rho_{\varepsilon}(\Gamma_{\widehat{\varepsilon_i}}) +  \sum_{i \in \Delta_4} (g_{\hat\imath}-1) +1.$$ 
\end{lemma}   
\dimo
Theorem 22 of \cite{Casali-Cristofori-Dartois-Grasselli} gives the following formula for the G-degree of $\Gamma$:
$$ \omega_G(\Gamma) \ = \ 6 \Big((p-1) - \sum_{i\in \Delta_4} (g_{\hat\imath} - 1) + (\chi(K(\Gamma))-2)\Big). $$
Hence the result comes by comparing it with formula \eqref{Euler-caracteristic n=4}. 
\qed

 \begin{lemma} \label{bound G-degree-order (n=4)}  
 Let $(\Gamma,\gamma)$ be an order $2p$ 5-colored graph representing a  compact $4$-manifold $M^4$ and satisfying $g_{\hat\imath}=1$ for each $i \in \Delta_4.$ 
 Then:   $$p =  1+ \frac{\omega_G(\Gamma) - \sum_{i \in \Delta_4}  \omega_G(\Gamma_{\hat\imath})} 3.$$
 In particular, 
 if $M^4$ is the product $\bar M \times I$, $\bar M$ being a closed 3-manifold: 
$$ \omega_G(\Gamma) \ge 5 k(\bar M).$$  
 \end{lemma}

 \dimo
 Since $g_{\hat\imath}=1$ $\forall i \in \Delta_4$,  \cite[Lemma 13]{Casali-Cristofori-Dartois-Grasselli}  directly yields 
 $$  \omega_G(\Gamma) = 3(p-1) + \sum_{i \in \Delta_4}  \omega_G(\Gamma_{\hat\imath}),$$
which proves the general statement. 

On the other hand, If $|K(\Gamma)|=\widehat M^4$ contains $h$  $(1 \le h \le 5)$ singular vertices, labelled with $h$ different colors, 
$\omega_G(\Gamma_{\hat\jmath_s}) \ge 3$ holds for (at least) $h$ colors $j_1, \dots, j_h \in \Delta_4.$   Hence,  $$ p \le \frac {\omega_G(\Gamma)} 3 - (h-1)$$ easily follows. 

In particular,  if $\Gamma$ represents $\bar M \times I$, the hypothesis $g_{\hat\imath}=1$ ($\forall i \in \Delta_4$) implies the existence of two colors $c_1, c_2 \in \Delta_4$ so that both the 4-residues $\Gamma_{\widehat {c_1}}$ and $\Gamma_{\widehat {c_2}}$ of $\Gamma$ represent $\bar M$. Then,  for each $i\in \{1,2\}$, \ $\omega_G(\Gamma_{\widehat {c_i}}) \ge \mathcal D_G(\bar M) = k(\bar M) = \bar p -1,$ with $\# V(\Gamma)= \# V(\Gamma_{\widehat {c_i}}) =2p \ge 2 \bar p$.  
Hence,  
$$  \omega_G(\Gamma) \ge 3(p-1) +  2 \mathcal D_G(\bar M) \ge 5(\bar p -1) = 5 \mathcal D_G(\bar M) = 5 k(\bar M). $$
\qed

 \bigskip

As a consequence, we are now able to classify all compact  $4$-manifolds up to G-degree  18 and - under a suitable condition - up to G-degree 
24.

\begin{prop} \label{classif_or n=4}
\par \noindent
Let $(\Gamma,\gamma)$ be a $5$-colored graph representing a compact $4$-manifold $M^4$. 
\begin{itemize}
\item[a.] If $ \omega_G(\Gamma) \in \{0,6\}$, \  then  \ $M^4 \cong \mathbb S^4;$
\item[b.]   if  $ \omega_G(\Gamma) = 12$, \  then  
\begin{itemize}  
\item[$\bullet$] either \ $ M^4 \in \{\mathbb S^4, \ \mathbb S^1 \times \mathbb S^3, \ \mathbb S^1 \tilde \times \mathbb S^3\}$
\item[$\bullet$] or \ $M^4 \in \{   \mathbb Y^4_1, \ \tilde{\mathbb Y}^4_1\};$
\end{itemize}
\item[c.]   if  $ \omega_G(\Gamma) = 18$, \ then   
\begin{itemize}  
\item[$\bullet$] either \  $M^4 \in \{\mathbb S^4, \ \mathbb S^1 \times \mathbb S^3, \ \mathbb S^1 \tilde \times \mathbb S^3\}$
\item[$\bullet$] or \  $M^4 \in \{ \mathbb Y^4_1, \ \tilde{\mathbb Y}^4_1\}$
\item[$\bullet$] or \  $M^4 \in \{L(2,1) \times I, \  (\mathbb S^1 \times \mathbb S^2) \times I,   \  (\mathbb S^1 \tilde \times \mathbb S^2) \times I \}.$
\end{itemize}
\end{itemize}
\noindent
No other 5-colored graph representing a compact $4$-manifold exists with $ \omega_G(\Gamma)          
\leq 23.$              
\par \noindent 
Moreover, if $(\Gamma,\gamma)$  has one singular color at most: 
\begin{itemize}
\item[d.] if   $ \omega_G(\Gamma) = 24$, \ then  
\begin{itemize}  
\item[$\bullet$] either \ $ M^4 \in \{\mathbb S^4, \mathbb S^1 \times \mathbb S^3, \mathbb S^1 \widetilde \times \mathbb S^3, \#_2(\mathbb S^1 \times \mathbb S^3),  \#_2(\mathbb S^1 \tilde \times \mathbb S^3),  \mathbb{CP}^2 \}$
\item[$\bullet$] or \  $M^4 \in \{ \mathbb Y^4_1, \tilde{\mathbb Y}^4_1,  \mathbb Y^4_1 \#  \mathbb Y^4_1, \mathbb Y^4_1 \#  \tilde{ \mathbb Y}^4_1,  \tilde{ \mathbb Y}^4_1 \#  \tilde{ \mathbb Y}^4_1,  
 {\mathbb Y}^4_2, \tilde{\mathbb Y}^4_2,$ 
 $(\mathbb S^1 \times \mathbb S^3) \#  \mathbb Y^4_1,$  $(\mathbb S^1 \tilde \times \mathbb S^3) \#  {\mathbb Y}^4_1,$ $(\mathbb S^1 \tilde \times \mathbb S^3) \#  \tilde {\mathbb Y}^4_1,$  
$\mathbb S^2 \times \mathbb D^2,$ $\xi_2\}. $
\end{itemize}
\end{itemize}
\end{prop}

\dimo If $\Gamma$ is a $5$-colored graph representing a compact $4$-manifold, then $ \omega_G(\Gamma) \le 6$ (resp.  $ \omega_G(\Gamma) \le 18$) 
(resp.  $ \omega_G(\Gamma) \le 24$) implies, via formula \eqref{multiplo6} (i.e.:  $  \omega_G(\Gamma) = 6 (\rho_\varepsilon (\Gamma) + \rho_{\varepsilon^{\prime}}(\Gamma))$, \ $(\varepsilon, \varepsilon^{\prime})$ being an arbitrary pair of associated permutations of $\Delta_4$),   the existence of a permutation $\varepsilon$ of $\Delta_4$ such that  $\rho_\varepsilon (\Gamma)=0$ 
(resp. $\rho_\varepsilon (\Gamma) \le 1$) (resp. $\rho_\varepsilon (\Gamma) \le 2$).  

In case   $\rho(\Gamma)= min_{\varepsilon} \ \rho_\varepsilon (\Gamma) =0$, \ $M^4 \cong \mathbb S^4$ follows from Proposition \ref{regular genus vs generalized regular genus} (b). 
This proves statement (a).  

\smallskip

In case   $\rho(\Gamma)= min_{\varepsilon} \ \rho_\varepsilon (\Gamma) =1$, Proposition \ref{gen-genus_one} ensures $\Gamma$ to represent either a closed $4$-manifold $M^4$, with $M^4 \in \{\mathbb S^4, \ \mathbb S^1 \times \mathbb S^3, \ \mathbb S^1 \tilde \times \mathbb S^3\}$,   or a compact $4$-manifold with non-empty boundary $M^4,$ with  $M^4 \in \{\mathbb Y^4_1, \ \tilde{\mathbb Y}^4_1, \bar M \times I\}$, \  $\bar M$ being a closed genus one 3-manifold. 
On the other hand, $\rho_\varepsilon (\Gamma) \le 1$  implies, via inequality \eqref{chiavacci-pareschi}, $\rho_{\varepsilon}(\Gamma_{\widehat{\varepsilon_i}}) \le \rho_\varepsilon (\Gamma) \le 1$ for each $i \in \Delta_4$;  hence, a (possible) sequence of proper 1-dipoles allows to consider the additional assumption $g_{\hat\imath}=1$ $\forall i \in \Delta_4$ (see Remark \ref{rem_dipoli}). Now, Lemma  \ref{bound G-degree-order (n=4)} ensures that,  if $M^4 \cong \bar M \times I,$ \ $\bar M$  ($\ne \mathbb S^4$) being a closed 3-manifold, then 
$ \omega_G(\Gamma) \ge 15$ holds.  
This fact completes the proof of statement (b). 

Moreover, Lemma  \ref{bound G-degree-order (n=4)}  proves that, if $M^4 \cong \bar M \times I$ and $ \omega_G(\Gamma) = 18$,  
the closed 3-manifold $\bar M$ has gem-complexity $k(\bar M) \le 3$.  The  existing classification of closed 3-manifolds via gem-complexity (see \cite{Casali-Cristofori JKTR 2008}) implies $\bar M$ to be either $L(2,1)$ or $\mathbb S^1 \times \mathbb S^2$ or $ \mathbb S^1 \tilde \times \mathbb S^2$. This completes the proof of statement  (c). 

\medskip

Let us now take into account the last case   $ \omega_G(\Gamma) = 24$, with $\rho_{\varepsilon}(\Gamma)= 2$ for each permutation $\varepsilon$ of $\Delta_4,$ and with the additional hypothesis that $\Gamma$ has one singular color at most. 
Some subcases occur: 
\begin{itemize}
\item{} $\exists c \in \Delta_4$ and a cyclic permutation $\varepsilon$ of $\Delta_4$ such that $\rho_{\varepsilon}(\Gamma_{\hat c}) =2.$   

In this case,  the identification of the represented compact $4$-manifold follows from $\rho_{\e}(\Gamma)=  \rho_{\varepsilon}(\Gamma_{\hat c}) =2 $, via Proposition \ref{null-subgenus}(b) 
 and  Proposition  \ref{gen-genus_two_onesingularcolor}: either   $M^4 \in \{\mathbb S^4, \ \mathbb S^1 \times \mathbb S^3, \ \mathbb S^1 \tilde \times \mathbb S^3,  \  \#_2(\mathbb S^1 \times \mathbb S^3), \  \#_2(\mathbb S^1 \tilde \times \mathbb S^3)\}$ (in case $K(\Gamma)$ has no singular vertex), or  $M^4 \in \{\mathbb Y^4_1, \tilde{\mathbb Y}^4_1,$ 
$(\mathbb S^1 \times \mathbb S^3) \#  \mathbb Y^4_1,$  $(\mathbb S^1 \tilde \times \mathbb S^3) \#  {\mathbb Y}^4_1,$ $(\mathbb S^1 \tilde \times \mathbb S^3) \#  \tilde {\mathbb Y}^4_1,$  
 ${\mathbb Y}^4_2, \tilde{\mathbb Y}^4_2\}$  (in case $K(\Gamma)$ has exactly one singular vertex),   or  $M^4 \in \{ \mathbb Y^4_1 \#  \mathbb Y^4_1,  \  \mathbb Y^4_1 \#  \tilde{ \mathbb Y}^4_1,  \  \tilde{\mathbb Y}^4_1 \#   \tilde{\mathbb Y}^4_1\}$  (in case $K(\Gamma)$ has exactly two singular vertices, labelled by the same color). 

\item{} For each permutation $\varepsilon$ of $\Delta_4,$ \ $\rho_{\varepsilon}(\Gamma_{\widehat{\varepsilon_i}}) \leq 1 \ \ \forall i\in\Delta_4$, and, therefore - via a (possible) sequence of proper 1-dipoles -, we may assume  $g_{\hat\imath}=1$ for each $i\in\Delta_4.$

If $M^4$ is closed, the result comes from Propositions 33 and 35 of \cite{Casali-Cristofori-Dartois-Grasselli}.

Otherwise, note that, by Lemma \ref{sommarho_i}, $4\leq p\leq 8$; hence the disjoint links of the singular vertices of $K(\Gamma)$ can only represent $L(2,1),\ \mathbb S^1\times \mathbb S^2$
or $\ \mathbb S^1\tilde\times \mathbb S^2$, since all other $4$-colored graphs with $p\leq 8$, not representing $\mathbb S^3$, admit at least one regular embedding into a surface of genus greater then one.
As a consequence, since $K(\Gamma)$ has exactly one singular vertex (in virtue of the assumption $\rho_{\varepsilon}(\Gamma_{\widehat{\varepsilon_i}}) \leq 1 \ \ \forall i\in\Delta_4$), Proposition \ref{gen-genus_two}  yields $M^4\in \{ \mathbb Y^4_1, 
\ \tilde{\mathbb Y}^4_1,  { \mathbb Y}^4_2, \ \tilde{ \mathbb Y}^4_2,$
$(\mathbb S^1 \times \mathbb S^3) \#  \mathbb Y^4_1,$  $(\mathbb S^1 \tilde \times \mathbb S^3) \#  {\mathbb Y}^4_1,$ $(\mathbb S^1 \tilde \times \mathbb S^3) \#  \tilde {\mathbb Y}^4_1,$    
$ \mathbb S^2 \times \mathbb D^2, \ \xi_2\}. $
\end{itemize}
\ \qed

\begin{rem} \label{1/N-expansion}
{\em As hinted to in Section 1, the G-degree $\omega_G$ of colored graphs plays an important role within colored tensor models theory, in virtue of the following $1/N$ expansion of the  correlation functions  
\begin{equation} \label{1/N expansion}
\frac{1}{N^d}\log \mathcal{Z}[N,\{t_{B}\}] = \sum_{\omega_G\ge 0}N^{-\frac{2}{(d-1)!}\omega_G}F_{\omega_G}[\{t_B\}]\in \mathbb{C}[[N^{-1}, \{t_{B}\}]],
\end{equation}
\noindent where the coefficients $F_{\omega_G}[\{t_B\}]$ are generating functions of connected bipartite $(d+1)$-colored graphs with fixed G-degree $\omega_G$. 
Moreover, \cite[Theorem 1]{Casali-Grasselli 2017} proves that, if $d$ even, $d\ge 4$,  the only non-null terms in the above formula are the ones corresponding to even (integer) powers of $1/N.$ 
\\
Hence, for $d=4$, Proposition \ref{classif_or n=4} yields the identification of all  compact orientable PL $4$-manifolds (resp. compact orientable  PL $4$-manifolds with empty or connected boundary), 
represented by regular graphs involved in the first four (resp. five) most significant non-null terms of the $1/N$ expansion of formula \eqref{1/N expansion}.}
\end{rem}

\medskip 

Proposition \ref{classif_or n=4} allows us to prove the second main result of the present paper (already stated in Section \ref{intro}).

\medskip

\dimo  (Theorem \ref{thm.classif_G-degree})\ \  It is well-known that the minimal (order ten) crystallizations of $\mathbb S^1 \times \mathbb S^3$ and $\mathbb S^1 \tilde \times \mathbb S^3$ have G-degree 12 (see \cite[Corollary 26]{Casali-Cristofori-Dartois-Grasselli}). 
Moreover,  in Proposition \ref{prop.handlebodies} $\mathbb Y^4_1$ and $\tilde{\mathbb Y}^4_1$ are proved to have G-degree 12, too.

Then it is easy to check that the 5-colored graphs obtained from the minimal (order eight) crystallizations of $L(2,1)$, $\mathbb S^1 \times \mathbb S^2$ and $\mathbb S^1 \tilde \times \mathbb S^2$ by applying the procedure described in Proposition \ref{MxI}
have G-degree equal to 18 and represent  $L(2,1) \times I$, $(\mathbb S^1 \times \mathbb S^2) \times I$ and $(\mathbb S^1 \tilde \times \mathbb S^2) \times I$ respectively. 

Further, the well-known 5-colored graph representing the closed $4$-manifolds $\#_2(\mathbb S^1 \times \mathbb S^3)$, $\#_2(\mathbb S^1 \tilde \times \mathbb S^3)$ and $\mathbb{CP}^2$ have G-degree $24$, as well as the 5-colored graphs representing ${\mathbb Y}^4_2$, $\tilde{ \mathbb Y}^4_2$,    $(\mathbb S^1 \times \mathbb S^3) \#  \mathbb Y^4_1,$  $(\mathbb S^1 \tilde \times \mathbb S^3) \#  {\mathbb Y}^4_1$ and $(\mathbb S^1 \tilde \times \mathbb S^3) \#  \tilde {\mathbb Y}^4_1$ obtained by graph connected sums  (see Section \ref{sec_general-properties}).

Finally, in Section \ref{section_bundles}, we have obtained a 5-colored graph $(\Lambda_0, \lambda_0)$ (resp. $(\Lambda_2, \lambda_2)$) with $\rho_\varepsilon(\Lambda_0)=\rho_{\varepsilon^{\prime}}(\Lambda_0) =  2$ (resp. $\rho_\varepsilon(\Lambda_2)=\rho_{\varepsilon^{\prime}}(\Lambda_2)=2$), where $(\varepsilon, \varepsilon^{\prime})$ is a pair of associated permutations of $\Delta_4$: see Figure \ref{fig.csi_c} (resp. Figure \ref{fig.S2xD2}). Hence, by  formula \eqref{multiplo6},              
$\omega_G(\Lambda_0)= \omega_G(\Lambda_2)= 6(2+2)=24$ follows. 

The statement is now a direct consequence of Proposition \ref{classif_or n=4}. 
\qed 

\begin{cor} \label{G-degree D2-fibrati}
Let  $\xi_c$ be the $\mathbb D^2$-bundle over $\mathbb S^2$ with Euler class $c$, $\forall c \in \mathbb  Z^+ - \{1\}$.  
Then: $$ { \mathcal D}_G(\xi_2)= {\mathcal D}_G(\mathbb S^2 \times \mathbb D^2)=24,$$
while 
$$ 30 \le  { \mathcal D}_G(\xi_c) \le 12 c \ \ \ \forall c \in \mathbb  Z^+ - \{1,2\}.$$  
Moreover: 
$$\mathcal D_G(\mathbb Y^4_1 \# \mathbb Y^4_1)  =  \mathcal D_G (\mathbb Y^4_1 \# \tilde{ \mathbb Y}^4_1)  = \mathcal D_G (\tilde{ \mathbb Y}^4_1 \#  \tilde{ \mathbb Y}^4_1) =  24.$$
\end{cor} 

\dimo
The statements concerning $\xi_c$,  $\forall c \in \mathbb  Z^+ - \{1\}$, and $\mathbb S^2 \times \mathbb D^2$ (resp. $\mathbb Y^4_1 \# \mathbb Y^4_1$, $\mathbb Y^4_1 \# \tilde{ \mathbb Y}^4_1$ and $\tilde{ \mathbb Y}^4_1 \#  \tilde{ \mathbb Y}^4_1$) are trivial consequence of  Theorem  \ref{thm.classif_G-degree}, together with the constructions presented in Section \ref{section_bundles} (resp. in Section \ref{sec_hadlebodies and xI}). 

Alternatively, ${ \mathcal D}_G(\xi_2)= {\mathcal D}_G(\mathbb S^2 \times \mathbb D^2)= {\mathcal D}_G(\mathbb Y^4_1 \# \mathbb Y^4_1)  =  {\mathcal D}_G (\mathbb Y^4_1 \# \tilde{ \mathbb Y}^4_1) =  {\mathcal D}_G (\tilde{ \mathbb Y}^4_1 \#  \tilde{ \mathbb Y}^4_1) = 24$ could also be proved directly from the computation of their generalized regular genus, performed in Corollary  \ref{D2-fibrati}, by making use of  the relation $ \mathcal D_G(N) \ \ge \ \frac{d!}2 \cdot \overline{\mathcal G}(N)$ (recalled in Section 2). 
\qed

\begin{conj} \label{congettura G-degree D2-fibrati}
Let  $\xi_c$ be the $\mathbb D^2$-bundle over $\mathbb S^2$ with Euler class $c$. 
Then: 
$$ { \mathcal D}_G(\xi_c) \ = \ 12 c \ \ \ \forall c \in \mathbb  Z^+ - \{1,2\}.$$
\end{conj}

The general computation of the G-degree of the products with the interval, performed in Section 4, gives rise in dimension 4 to the following result. 

\begin{prop}  \label{G-degree MxI (n=4)}
For each closed 3-manifold $M$,   $$\mathcal D_G (M \times I) \le 6 \cdot \mathcal D_G (M).$$

\end{prop}

\dimo   
If $(\Gamma, \gamma)$ is a crystallization of $M$ realizing gem-complexity (i.e.: $\#V(\Gamma)=2 \bar p$, where $k(M) = \bar p -1$), it is well-known  that $\mathcal D_G (M) = \omega_G (\Gamma)= \bar p -1$ (\cite{Casali-Cristofori-Dartois-Grasselli}).   
If $\tilde \Gamma$ is the $5$-colored graph representing $M \times I$ considered in Proposition \ref{MxI}, the last formula of Proposition \ref{MxI} becomes: 
$$\omega_G(\tilde \Gamma)  \ = \ 3  \Big[\sum_{i \in \{1, 2, 3\}} (\bar p - g_{\varepsilon_0 \varepsilon_i}) - (\bar p-1) + \omega_G(\Gamma)\Big], $$
$\varepsilon = (\varepsilon_0 ,\varepsilon_1, \varepsilon_2, \varepsilon_3)$ being the cyclic permutation of $\Delta_3$ so that $\rho(\Gamma)= \rho_{\varepsilon}(\Gamma)$. 
On the other hand, $\sum_{i \in \{1, 2, 3\}} g_{\varepsilon_0 \varepsilon_i} = g_{\varepsilon_0 \varepsilon_1} + g_{\varepsilon_0 \varepsilon_2} + g_{\varepsilon_0 \varepsilon_3} = g_{\varepsilon_0 \varepsilon_1} + g_{\varepsilon_0 \varepsilon_2} + g_{\varepsilon_1 \varepsilon_2} =  2 - 2 \rho(\Gamma_{\widehat{\varepsilon_3}})  + \bar p = 2 + \bar p$, 
since both $g_{i,j} = g_{\hat\imath \hat\jmath}$ and $\rho(\Gamma_{\hat\imath})=0$ ($\forall i,j \in \Delta_3$) hold  in any crystallization of a closed 3-manifold.   

Hence:
$$\omega_G(\tilde \Gamma)  \ = \ 3  \Big[3\bar p - (2 + \bar p)  - (\bar p-1) + \omega_G(\Gamma)\Big] = 6 \cdot (\bar p -1) = 6 \cdot \mathcal D_G (M).$$
The claim now trivially follows.
\qed

\begin{cor} \label{G-degree L(3,1)xI}
$$ { \mathcal D}_G(L(3,1) \times I) = 30.$$
\end{cor} 

\dimo
Since  $\mathcal D_G (L(3,1))= k(L(3,1))=5$ (see \cite{Casali-Cristofori JKTR 2008} and \cite[Theorem 16]{Casali-Cristofori-Dartois-Grasselli}) and $\mathcal D_G(N) \equiv 0 \mod 6$ for each singular 4-manifold (see \cite{Casali-Grasselli 2017} or formula \eqref{multiplo6}),  
the statement is a trivial consequence of  Proposition \ref{G-degree MxI (n=4)}, together with Theorem  \ref{thm.classif_G-degree} (or Lemma \ref{bound G-degree-order (n=4)}).   
\qed

\noindent\textbf{Acknowledgments.}  This work was supported by GNSAGA of INDAM and by the University of Modena and Reggio Emilia, project: {\it ``Discrete Methods in Combinatorial Geometry and Geometric Topology".}


\begin{thebibliography}{9}
\bibitem{BCrG}
P.Bandieri, P. Cristofori, and C.Gagliardi, Nonorientable 3-manifolds admitting coloured triangulations with at most 30 tetrahedra. \textit{J. Knot Theory Ramifications}~{\bf 18} (2009), 381-395.

\bibitem{Basak-Casali 2016}
B.Basak, and M.R.Casali, Lower bounds for regular genus and gem-complexity of PL 4-manifolds. \textit{Forum Math.}~{\bf 29} (2017), no. 4, 761-773. 

\bibitem{[uncoloring]}
V. Bonzom,  R. Gurau,  V. Rivasseau, Random tensor models in the large N limit: Uncoloring the colored tensor models. \textit{Phys. Rev. D}~{\bf 85} 084037 (2012).

\bibitem{[Bonzom-Gurau-Riello-Rivasseau]}
V. Bonzom,  R. Gurau,  A. Riello, and V. Rivasseau, Critical behavior of colored tensor models in the large N limit. \textit{Nucl. Phys. B}~{\bf 853} (2011), no. 1, 174-195.

\bibitem{[Bonzom-Lionni-Tanasa]} 
V. Bonzom,  L. Lionni, and A. Tanasa, Diagrammatics of a colored SYK model and of an SYK-like tensor model, leading and next-to-leading orders. \textit{J. Math. Phys.}~{\bf 58} 052301 (2017).

\bibitem{Carrozza-Tanasa}
S. Carrozza, and A. Tanasa, O(N) Random Tensor Models. \textit{Lett. Math. Phys.}~{\bf 106} (2016), no. 11, 1531-1559.

\bibitem{Casali - Forum Math 1992}
M. R. Casali, A combinatorial characterization of 4-dimensional handlebodies. \textit{Forum Math.}~{\bf 4} (1992), 123-134.

\bibitem{Casali UMI}
M. R. Casali, An infinite class of bounded 4-manifolds having regular genus three. \textit{Boll. Unione Mat. Ital.}~{\bf 10-A} (1996), 279-303.

\bibitem{Casali Canadian}
M. R. Casali, Classifying PL 5-manifolds by regular genus: the boundary case. \textit{Canad. J. Math.}~{\bf 49} (1997), 193-211. 

\bibitem{Casali_Forum2003} 
M. R. Casali, On the regular genus of 5-manifolds with free fundamental group. \textit{Forum Math.}~{\bf 15} (2003), 465-475.  

\bibitem{Casali-Cristofori JKTR 2008}
M. R. Casali, and P. Cristofori, A catalogue of orientable 3-manifolds triangulated by $30$ coloured tetrahedra. \textit{J. Knot Theory Ramifications}~{\bf 17} (2008), 1-23.

\bibitem{Casali-Cristofori ElecJComb 2015}
M. R. Casali, and P. Cristofori, Cataloguing PL 4-manifolds by gem-complexity. \textit{Electron J. Combin.}~{\bf 22} (2015), no. 4, \#P4.25.

\bibitem{[Casali-Cristofori_trisection]}
M. R. Casali, and P. Cristofori, Gem-induced trisections of compact PL 4- manifolds, 2022.  arXiv:1910.08777v3
 
\bibitem{Casali-Cristofori-Gagliardi Complutense 2015}
M. R. Casali,  P. Cristofori, and C. Gagliardi, {\it Classifying PL 4-manifolds via crystallizations: results and open problems}. 
In  {\it Mathematical Tribute to Professor Jos\'e Mar\'ia Montesinos Amilibia}, Universidad Complutense Madrid, 2016. \ [ISBN: 978-84-608-1684-3]

\bibitem{Casali-Cristofori-Dartois-Grasselli}
M. R. Casali,  P. Cristofori,  S. Dartois, and L. Grasselli, Topology in colored tensor models  via crystallization theory. \textit{J. Geom. Phys.}~{\bf 129} (2018), 142-167. 

\bibitem{Casali-Cristofori-Grasselli}
M. R. Casali,  P. Cristofori, and L. Grasselli, G-degree for singular manifolds. \textit{RACSAM}~{\bf 112} (2018), no. 3, 693-704. 

\bibitem{Casali-Gagliardi ProcAMS}
M. R. Casali, and  C. Gagliardi, Classifying PL 5-manifolds up to regular genus seven. \textit{Proc. Amer. Math. Soc.}~{\bf 120}  (1994),  275-283.

\bibitem{Casali-Grasselli 2017}
M. R. Casali, and L. Grasselli, Combinatorial properties of the G-degree. \textit{Rev. Mat. Complut.}~{\bf 32} (2019), no. 1, 239-254.

\bibitem{Casali-Malagoli}
M. R. Casali, and L. Malagoli, Handle-decompositions of PL 4-manifolds. \textit{Cah. Topol. G\'eom. Diff\'er. Cat\'eg.}~{\bf 38}  (1997), 141-160.

\bibitem{Chiavacci}
R. Chiavacci, Pseudocomplessi colorati e loro gruppi fondamentali. \textit{Ric. Mat.}~{\bf  35} (1986), 247-268.

\bibitem{Chiavacci-Pareschi}
R. Chiavacci, and G. Pareschi, Some bounds for the regular genus of closed PL manifolds. \textit{Discrete Math.}~{\bf 82} (1990), 165-180.

\bibitem{Cristofori-Fomynikh-Mulazzani-Tarkaev}
P. Cristofori,  E. Fomynikh,  M. Mulazzani, and V. Tarkaev, 4-colored graphs and knot/link complements. \textit{Results Math.}~{\bf 72} (2017), no. 1-2, 471-490.

\bibitem{Cristofori-Gagliardi-Grasselli}
P. Cristofori,  C. Gagliardi, and L. Grasselli, Heegaard and regular genus of 3-manifolds with boundary. \textit{Rev. Mat. Complut.}~{\bf 8}  (1995), 379-398. 

\bibitem{Cristofori-Mulazzani}
P. Cristofori, and M. Mulazzani, Compact 3-manifolds via 4-colored graphs. \textit{RACSAM}~{\bf 110} (2015), no. 2, 395-416.

\bibitem{DiFrancesco-al}
P. Di Francesco, P. Ginsparg, and J. Zinn-Justin, 2D gravity and random matrices. \textit{Phys. Rep.}~{\bf 254} (1995), 1-133.

\bibitem{Ferri-Gagliardi Proc AMS 1982}
M. Ferri, and C. Gagliardi, The only genus zero n-manifold is $\mathbb S^n$. \textit{Proc. Amer. Math. Soc.}~{\bf 85} (1982), 638-642.

\bibitem{Ferri-Gagliardi Yokohama 1985}
M. Ferri, and C. Gagliardi, A characterization of punctured n-spheres. \textit{Yokohama Math. J.}~{\bf 33} (1985), 29-38.

\bibitem{Ferri-Gagliardi-Grasselli} 
M. Ferri,  C. Gagliardi, and L. Grasselli, A graph-theoretical representation of PL-manifolds. A survey on crystallizations. \textit{Aequationes Math.}~{\bf 31} (1986), 121-141.

\bibitem{Fusy-Lionni-Tanasa}
E. Fusy, L. Lionni, and A. Tanasa, Combinatorial study of graphs arising from the Sachdev-Ye-Kitaev model. \textit{European J. Combin.}~{\bf 86} (2020), 103066.

\bibitem{Gabai}
D. Gabai, Foliations and the topology of 3-manifolds. III. \textit{J. Differential Geom.}~{\bf 26} (1987), 479-536.
 
 \bibitem{Gagliardi 1981}
C. Gagliardi, Extending the concept of genus to dimension $n$. \textit{Proc. Amer. Math. Soc.}~{\bf 81} (1981), 473-481.

\bibitem{Gagliardi 1987}
C. Gagliardi, On a class of 3-dimensional polyhedra. \textit{Ann. Univ. Ferrara}~{\bf 33} (1987), 51-88.

\bibitem{Gagliardi_boundary}
C. Gagliardi, Regular genus: the boundary case. \textit{Geom. Dedicata}~{\bf 22} (1987), 261-281. 

\bibitem{Gordon-Luecke}
C. McA. Gordon, and J. Luecke, Knots are determined by their complements. \textit{Bull. Amer. Math. Soc.}~{\bf 20} (1989), no. 1, 83-87.

\bibitem{Grasselli-Mulazzani}
L. Grasselli, and M. Mulazzani, Compact $n$-manifolds via $(n+1)$-colored graphs: a new approach. \textit{Algebra Colloq.}~{\bf 27} (2020), no. 1, 95-120.   

\bibitem{Gurau-book}
R. Gurau, {\it Random Tensors}, Oxford University Press, 2016.

\bibitem{Gurau2017}
R. Gurau, The complete $1/N$ expansion of a SYK-like tensor model. \textit{Nucl. Phys. B}~{\bf 916} (2017), 386-401. 

\bibitem{Gurau2019}
R. Gurau, Notes on Tensor Models and Tensor Field Theories. {\it Ann. Inst. Henri Poincar\'e D} {\bf 9} (2022),  no. 1, 159-218.

\bibitem{Gurau-Ryan}
R. Gurau, and J. P. Ryan, Colored Tensor Models - a review. \textit{SIGMA}~{\bf 8} (2012),  020.

\bibitem{[Gurau-Schaeffer 2013]}
R. Gurau, and G. Schaeffer, Regular colored graphs of positive degree. \textit{Ann. Inst. Henri Poincar\'e D}~{\bf 3} (2016), 257-320.

\bibitem{Koda-Martelli-Naoe} 
Y. Koda,  H. Naoe, and B. Martelli, Four-manifolds with shadow-complexity one.  {\it Ann. Fac. Sci. Toulouse} {\bf 31}  (2022), no. 4, 1111-1212. 

\bibitem{Kronheimer}
P. Kronheimer,  T. Mrowka,  P. Ozsv\'ath, and Z. Szab\'o, Monopoles and Lens Space Surgeries. \textit{Annals of Math.}~{\bf 165} (2007), no. 2, 457-546. 

\bibitem{Laudenbach-Poenaru}
F. Laudenbach, and V. Poenaru, A note on 4-dimensional handlebodies. \textit{Bull. Soc. Math. France}~{\bf 100} (1972), 337-344.

\bibitem{Martelli} 
B. Martelli, Four-manifolds with shadow-complexity zero. \textit{Int. Math. Res. Not.}~{\bf 2011} (2011), no. 6, 1268-1351.  

\bibitem{Meier} 
J. Meier, Trisections and spun 4-manifolds. \textit{Math. Res. Lett.}~{\bf 25} (2018), no. 5, 1497-1524. 

\bibitem{Meier-Zupan} 
J. Meier, and A. Zupan, Genus-two trisections are standard. \textit{Geom. Topol.}~{\bf 21} (2017), no. 3, 1583-1630. 

\bibitem{Montesinos} 
J. M. Montesinos Amilibia, {\it Heegaard diagrams for closed 4-manifolds}.  In {\it Geometric topology}, Proc. 1977 Georgia Conference, Academic Press, 1979, 219-237. 

\bibitem{Naoe-Proc.AMS} 
H. Naoe, Shadows of 4-manifolds with complexity zero and polyhedral collapsing. \textit{Proc. Amer. Math. Soc.}~{\bf 145}  (2017), 4561-4572. 

\bibitem{Pezzana}
M. Pezzana, Sulla struttura topologica delle variet\`a compatte. \textit{Atti Semin. Mat. Fis. Univ. Modena}~{\bf 23} (1974), 269-277.

\bibitem{Spreer-Tillmann} 
J. Spreer, and S. Tillmann, The trisection genus of standard simply connected PL 4-manifolds. 
\textit{34nd International Symposium on Computational Geometry (SoCG 2018)}. In Leibniz International Proceedings in Informatics (LIPICS), 99, 2018, 71:1-71:13.  

\bibitem{Tanasa-multiorientable}
A. Tanasa, The Multi-Orientable Random Tensor Model, a Review. \textit{SIGMA}~{\bf 12} (2016), 056, 23 pages.

\bibitem{Witten} E. Witten, An SYK-like model without disorder. \textit{Journal of Physics A: Mathematical and Theoretical}~{\bf 52} (47) 474002 (2019). 


\end{thebibliography}
\end{document}